\documentclass[11pt]{article}

\usepackage{graphicx}
\usepackage{color}
\usepackage{amsmath,amssymb,amsthm,amsfonts}
\usepackage{amssymb}
\newtheorem{thm}{Theorem}[section]

\newtheorem{lem}[thm]{Lemma}
\newtheorem{prop}[thm]{Proposition}

\theoremstyle{definition}

\theoremstyle{remark}

\numberwithin{equation}{section}

\DeclareMathSymbol{\C}{\mathalpha}{AMSb}{"43}

\newcommand{\eps}{\varepsilon}

\newcommand{\alp}{\alpha}

\newcommand{\R}{{\mathbb{R}}}

\newcommand{\bsub}{\begin{subequations}}
\newcommand{\esub}{\end{subequations}$\!$}

\begin{document}
\title{Minimizers of $L^{2}$-critical inhomogeneous  variational problems with a spatially decaying nonlinearity in bounded domains}


\author{Hongfei Zhang\thanks{School of Mathematics and Statistics, Central China Normal University, P.O. Box 71010, Wuhan 430079,
P. R. China.  Email: \texttt{hfzhang@mails.ccnu.edu.cn}.}
 \, and\, Shu Zhang\thanks{School of Mathematics and Statistics, Central China Normal University, P.O. Box 71010, Wuhan 430079,
P. R. China.  Email: \texttt{shu@mails.ccnu.edu.cn}.}
}

\date{\today}

\smallbreak \maketitle

\begin{abstract}

We consider the minimizers of $L^{2}$-critical inhomogeneous variational problems with a spatially decaying nonlinear term in an open bounded domain $\Omega$ of $\mathbb{R}^{N}$ which contains $0$. We prove that there is a threshold $a^{*}>0$ such that minimizers exist for $0<a<a^{*}$ and the minimizer does not exist for any $a>a^{*}$. In contrast to the homogeneous case,  we show that both the existence and nonexistence of minimizers may occur at the threshold $a^*$ depending on the value of $V(0)$, where $V(x)$ denotes the trapping potential. Moreover, under some suitable assumptions on $V(x)$, based on a detailed analysis on the concentration behavior of minimizers as $a\nearrow a^*$, we prove local uniqueness of minimizers when $a$ is close enough to $a^*$.
\end{abstract}
\vskip 0.2truein

\section{Introduction}

In this paper, we consider the following $L^{2}$-critical constraint inhomogeneous variational problem
\begin{equation}\label{1.1}
e(a):=\inf_{u\in \mathcal{M}}E_{a}(u),
\end{equation}
where the energy functional $E_{a}(u)$ contains a spatially decaying nonlinearity and is defined by
\begin{equation}\label{1.2}
E_{a}(u):=\int_{\Omega}(|\nabla u|^{2}+V(x)|u|^{2})dx-\frac{a}{1+\beta^{2}}\int_{\Omega}\frac{|u|^{2+2\beta^{2}}}{|x|^{b}}dx,
\end{equation}
and the space $\mathcal{M}$ is defined as
\begin{equation}\label{1.3}
\mathcal{M}:=\{u\in H_{0}^{1}(\Omega):\int_{\Omega}|u(x)|^{2}dx=1\}.
\end{equation}
Here we assume that $N\geq1$, $a>0$, $0<b<min\{2,N\}$ and $\beta=\sqrt{\frac{2-b}{N}}$ is the $L^{2}$-critical exponent. Moreover, we suppose that the bounded domain $\Omega$ of $\mathbb{R}^{N}$ and the trapping potential $V(x)\geq0$ satisfy the following assumptions respectively.
\begin{itemize}
\item[\rm($V_{0})$.] $\Omega$ is an open connected domain containing $0$ and $\partial\Omega\in C^{1}$.
\item[\rm($V_{1})$.] $V(x)\in C(\bar{\Omega})$ and $\min_{x\in\overline{\Omega}}V(x)=0$.
\end{itemize}
From now on, we assume that $\Omega$ satisfies $(V_{0})$ throughout this paper.

The variational problem (\ref{1.1}) arises in various physical contexts, including the propagation of a laser beam in the optical fiber, Bose-Einstein condensates (BECs), and nonlinear optics (cf.\cite{GPA, BP, LT}). When $b=0$, (\ref{1.1}) is a homogeneous constraint variational problem. For $b=0$ and $N=2$, Luo and Zhu proved that there exists a critical constant $a_{0}>0$ such that minimizers of $e(a)$ exist if and only if $0<a<a_{0}$, and the limit behavior of minimizers as $a\nearrow a_{0}$ is also analyzed in \cite{ZL}. When $b>0$, the variational problem (\ref{1.1}) contains the inhomogeneous nonlinear term $K(x)|u|^{2+2\beta^{2}}$, where $K(x)=\frac{1}{|x|^{b}}$ admits $x=0$ as a singular point.
Compared with the homogeneous case where $b=0$, we shall point out that our analysis on inhomogeneous problem (\ref{1.1}) is more complex due to the singularity of $|x|^{-b}$.

If the bounded domain $\Omega$ is replaced by the whole space $\mathbb{R}^{N}$, the existence and nonexistence of minimizers for (\ref{1.1}) and various quantitative properties of (\ref{1.1}) were investigated extensively, see
\cite{AA, DGL1, DGL2, AI, SC, GS, GWZ, GZZ, LZ} and the references therein. Note that most of the problems in these studies are homogeneous. Recently, the inhomogeneous $L^{2}$-constraint variation problems were also analyzed in \cite{DGL1, DGL2}, where the authors addressed the existence, nonexistence and limit behavior of minimizers mainly in the case that the inhomogeneous nonlinear term $f(x)|u(x)|^{2+2\beta^{2}}$ with $f(x)\in L^{\infty}(\mathbb{R}^{N})$.
More recently, Luo and Zhang in \cite{LZ} analyzed the variational problem (\ref{1.1}) defined in $\mathbb{R}^{N}$, where the authors made full use of the following nonlinear Schr\"{o}dinger Equation
\begin{equation}\label{1.5}
-\Delta u+u-|x|^{-b}|u|^{2\beta^{2}}u=0\ \text{ in}\,\ \ \mathbb{R}^N.
\end{equation}
Compared with \cite{LZ}, we emphasize that the bounded domain $\Omega$ is variable under the scaling transformation and the optimal constant of the Gagliardo-Nirenberg inequality for the bounded domain $\Omega$ is not attained. Therefore, we need to explore some new analytic approaches to overcome these difficulties.

It is well-known from \cite{F,G,TJ} that (\ref{1.5}) admits a unique positive radial solution for any $N\geq1$. Such a positive solution of (\ref{1.5}) is always denoted by $Q=Q(|x|)$ throughout this paper. Moreover, we cite from \cite{AI} the following Gagliardo-Nirenberg inequality
\begin{equation}\label{1.6}
\frac{a^{*}}{1+\beta^{2}}\int_{\mathbb{R}^{N}}|x|^{-b}|u(x)|^{2+2\beta^{2}}dx\leq
\int_{\mathbb{R}^{N}}|\nabla u(x)|^{2}dx\Big(\int_{\mathbb{R}^{N}}|u(x)|^{2}dx\Big)^{\beta^{2}}\ \text{ in}\,\ \ H^{1}(\mathbb{R}^N),
\end{equation}
where the constant $a^{*}>0$ is given by
\begin{equation}\label{1.7}
a^{*}:=\|Q\|_{L^{2}(\mathbb{R}^{N})}^{2\beta^{2}},
\end{equation}
and the identity of $(\ref{1.6})$ holds if and only if $u(x)=mn^{\frac{N}{2}}Q(nx)$ ($m, n\neq 0$ are arbitrary). Recall also from \cite[Theorem 2.2]{SC} that $Q$ admits the following exponential decay
\begin{equation}\label{1.8}
|Q(x)|,\ \ |\nabla Q(x)|\leq Ce^{-|x|}\ \ \text{as}\  \ |x|\rightarrow\infty.
\end{equation}
Moreover, we can derive from (\ref{1.5}) and (\ref{1.6}) that $Q$ satisfies
\begin{equation}\label{1.9}
\int_{\mathbb{R}^{N}}|\nabla Q|^{2}dx=\frac{1}{\beta^{2}}\int_{\mathbb{R}^{N}}|Q|^{2}dx
=\frac{1}{1+\beta^{2}}\int_{\mathbb{R}^{N}}|x|^{-b}|Q|^{2+2\beta^{2}}dx.
\end{equation}
We next define the linearized operator $\mathcal{L}$ by
\begin{equation}\label{A1}
\mathcal{L}:=-\Delta +1-(1+2\beta^2)|x|^{-b}Q^{2\beta^2}\ \  \mbox{in}\ \, \R^N,
\end{equation}
where $Q=Q(|x|)>0$ is the unique positive  solution of (\ref{1.5}). It then follows from \cite{SC} that
\begin{equation}\label{4.1}
ker (\mathcal{L})=span \{0\}\quad\hbox{if}\,\ N\geq 3.
\end{equation}

Stimulated by the works mentioned above, the main purpose of this paper is to
address the existence and nonexistence of minimizers for \eqref{1.1}, and analyze the local uniqueness of minimizers for a special class of $V (x)$.

Before stating our results, we first introduce the Gagliardo-Nirenberg type inequality: for any open bounded domain $\Omega\subset\mathbb{R}^{N}$, it holds
\begin{equation}\label{1.10}
\frac{a^{*}}{1+\beta^{2}}\int_{\Omega}|x|^{-b}|u(x)|^{2+2\beta^{2}}dx\leq\int_{\Omega}|\nabla u(x)|^{2}dx\Big(\int_{\Omega}|u(x)|^{2}dx\Big)^{\beta^{2}},\ \ u\in H_{0}^{1}(\Omega),
\end{equation}
where the optimal constant $\frac{1+\beta^{2}}{a^{*}}$ is not attained. The proof of (\ref{1.10}) is left to the Appendix. Applying the inequality \eqref{1.10}, we shall establish the following existence and nonexistence of minimizers for (\ref{1.1}) under the general assumption $(V_{1})$.

\begin{thm}\label{th1}
Let $N\geq1$, $0<b<min{\{2,N\}}$, $\beta^{2}=\frac{2-b}{N}$, and $V(x)$ satisfies $(V_{1})$, then
\begin{enumerate}
\item If $0<a<a^{*}:=\|Q\|_{L^{2}(\mathbb{R}^{N})}^{2\beta^{2}}$, there exists at least one minimizer for (\ref{1.1}).

\item If $a>a^{*}$, there is no minimizer for (\ref{1.1}).
\end{enumerate}
\end{thm}

It is worth noting that the threshold $a^{*}$ stated in Theorem \ref{th1} is independent of the external potential $V(x)$ and the domain $\Omega$ as well. Moreover, in view of Theorem \ref{th1}, it is therefore natural to wonder whether the minimizers of $e(a)$ exist for the case $a=a^{*}$. The following result presents that $e(a)$ may admit minimizers at the threshold $a^{*}$ for some trapping potentials $V(x)$ satisfying $(V_{1})$. Indeed, using the Gagliardo-Nirenberg inequality (\ref{1.10}), one can check that $e(a^{*})\geq0$. On the other hand, taking the same trial function as in (\ref{eq2.8}) below and letting $\tau\rightarrow\infty$, we can obtain that $e(a^{*})\leq V(0)$. Therefore, we conclude from above that
$$0\leq e(a^{*})\leq V(0),$$
and our second result can be summarized as follows.


\begin{thm}\label{th2}
Let $N\geq1$, $0<b<min{\{2,N\}}$, $\beta^{2}=\frac{2-b}{N}$, and $V(x)$ satisfies $(V_{1})$, then we have
\begin{enumerate}
\item If $0\leq e(a^{*})< V(0)$, then there exists at least one minimizer for $e(a^{*})$.
\item If $V(0)=0$, then there is no minimizer for $e(a^{*})$. Moreover, $lim_{a\nearrow a^{*}}e(a)=e(a^{*})=V(0)$.
\end{enumerate}
\end{thm}



In order to analyze the uniqueness of minimizers for $e(a)$ when $a$ is close enough to $a^*$, we first discuss the limit behavior of minimizers as $a\nearrow a^*$. Let $u_{a}$ be a minimizer of $e(a)$. It then follows from the variational theory that $u_{a}$ satisfies the following nonlinear Schr\"{o}dinger equation
\begin{equation}\label{1.4}
-\Delta u_{a}+V(x)u_{a}-a|x|^{-b}|u_{a}|^{2\beta^{2}}u_{a}=\mu_{a} u_{a}\ \text{ in}\,\ \ \Omega,
\end{equation}
where $\mu_{a}\in\mathbb{R}$ is a suitable Lagrange multiplier associated to $u_{a}$.

For any minimizer $u_a$ of $e(a)$, it yields from \cite[Theorem 6.17]{LE} that $E_{a}(u_a)=E_{a}(|u_a|)$, which further implies that $|u_a|$ is a minimizer of $e(a)$. By the strong maximum principle, we can derive from (\ref{1.4}) that $|u_a|>0$ holds in $\Omega$. Therefore, without loss of generality, we can restrict the minimizers of $e(a)$ to positive functions. Motivated by \cite{GL,ZL}, in order to analyze the limit behavior of $u_{a}$ as $a\nearrow a^*$, some additional assumptions on $V(x)$ are required.

We shall assume that the external potential $V(x)\in C(\bar{\Omega})$ has $n\geq1$ distinct zero points $x_{i}\in\bar{\Omega}$ satisfying $V(x_{i})=0$ and $V(x)>0$ for any $x\neq x_{i}$, where $i=0,1,\cdots,n-1$. Considering Theorem \ref{th2} (2), we suppose that $x_{0}=0$ is a minimum point of $V(x)$. More precisely, we suppose that there exist $p_{i}>0$ and $C>0$ such that
\begin{enumerate}
\item [\rm($V_2$).]
$V(x)=h(x)|x|^{p_{0}}\prod_{i=1}^{n-1}|x-x_{i}|^{p_{i}}$\ \ with\ \ $C<h(x)<\frac{1}{C}$\ \ for\ \ all\ \ $x\in\bar{\Omega}$,
and $\lim_{x\rightarrow x_{i}}h(x)$ exists for all $0\leq i\leq n-1$.
\end{enumerate}
 For convenience, we denote
\begin{equation}\label{1.11}
p:=p_{0}>0,
\end{equation}
and
\begin{equation}\label{1.12}
\lambda:=\Big(\frac{p}{2}\int_{\mathbb{R}^{N}}|x|^{p}|Q(x)|^{2}dx\lim_{x\rightarrow0}
(h(x)\prod_{i=1}^{n-1}|x-x_{i}|^{p_{i}})\Big)^{\frac{1}{p+2}}>0.
\end{equation}
By virtue of the conditions on $V(x)$, the main result of concentration behavior is stated as follows.
\begin{thm}\label{th3}
Let $N\geq1$, $0<b<min{\{2,N\}}$, $\beta^{2}=\frac{2-b}{N}$, $0<a<a^*$, and $V(x)$ satisfies $(V_{1})$ and $(V_{2})$. Then for any positive minimizer $u_{a}$ of (\ref{1.1}), we have
\begin{equation}\label{1.13}
\lim_{a\nearrow a^{*}}
\Big(\frac{(a^*-a)\|Q\|_{2}^{2}}{a^*\beta^{2}\lambda^{p+2}}\Big)^{\frac{N}{2(p+2)}}
u_{a}\Big(\Big(\frac{(a^*-a)\|Q\|_{2}^{2}}{a^*\beta^{2}}\Big)^{\frac{1}{p+2}}
\frac{x}{\lambda}
\Big)=\frac{Q(x)}{\|Q\|_{2}}
\end{equation}
strongly in $H^1(\R^N)\cap L^{\infty}(\R^{N})$, where $p>0$ and $\lambda>0$ are defined by (\ref{1.11}) and (\ref{1.12}) respectively.
\end{thm}

We remark that the limit (\ref{1.13}) is defined in $H^1(\R^N)\cap L^{\infty}(\R^{N})$ in the sense that $u_{a}(x)\equiv0$ for all $x\in\mathbb{R}^{N}\setminus\Omega$. Because the limit equation \eqref{1.5} has no translation invariance,  it is worth pointing out that $u_{a}$ must concentrate at the origin rather than any global  minimum point of $V(x)$ as $a\nearrow a^{*}$. Moreover, Theorem \ref{th3} gives explicitly concentration rates of $u_{a}$ as $a\nearrow a^{*}$, that is
\begin{equation*}
u_{a}(x)\approx\frac{1}{\|Q\|_{2}}\Big(\frac{a^*\beta^2\lambda^{p+2} }{(a^*-a)\|Q\|_{2}^{2}}\Big)^{\frac{N}{2(p+2)}}Q\Big(\Big(\frac{a^*\beta^2}
{(a^*-a)\|Q\|_{2}^{2}}\Big)^{\frac{1}{p+2}}\lambda x\Big)\ \ \text{ as } \ \  a\nearrow a^*,
\end{equation*}
where $f(a)\approx g(a)$ means that $f/g\rightarrow1$ as $a\nearrow a^*$.

Theorem \ref{th3} is proved by the blow up analysis, for which we need to derive the exact blow up rate of minimizers as $a\nearrow a^{*}$. We shall reach this purpose by analyzing the refined energy estimates of $e(a)$ as $a\nearrow a^{*}$. However, due to the singularity of $|x|^{-b}$ and the variability of the scaling transformation of the bounded domain $\Omega$, we remark that the standard elliptic regularity theory is used with caution when investigating the $L^{\infty}$-uniform convergence of (\ref{1.13}).


Motivated by \cite{YPC,WLG, LL} and the references therein, based on Theorem \ref{1.3}, we finally investigate the uniqueness of positive minimizers for $e(a)$ as $a\nearrow a^*$ under some suitable assumptions on $V(x)$. 
\begin{thm}\label{th4}
Let $N\geq3$, $0<b<min\{2,\frac{N}{2}\}$, $\beta^{2}=\frac{2-b}{N}$ and assume that $V(x)$ satisfies $(V_{1})$ and $(V_{2})$. Then there exists a unique positive minimizer of $e(a)$ as $a\nearrow a^*$.
\end{thm}

We note that if the non-degeneracy property (\ref{4.1}) still holds for all $N\geq1$, then the restriction $N\geq3$ in Theorem \ref{th4} can be removed. In addition, the condition $0<b<min\{2,\frac{N}{2}\}$ guarantees that the integral \eqref{4.28} makes sense. Inspired by \cite{YPC,WLG, LL}, we shall prove Theorem \ref{th4} by constructing Pohozaev identity, which is widely used in the existing literature of studying the real-valued elliptic PDEs. However, the calculations involved in the proof are more complicated due to the existence of inhomogeneous nonlinear term.


This paper is organized as follows. In Section 2 we shall prove Theorem \ref{th1} and Theorem \ref{th2} on the existence and nonexistence of minimizers for $e(a)$. In Section 3, we first establish some energy estimates of positive minimizers for $e(a)$ as $a\nearrow a^*$, after which we complete the proof of Theorem \ref{th3}. By establishing local Pohozaev identity, we then prove Theorem \ref{th4}  on the
local uniqueness of positive minimizers in Section 4.  In the Appendix, we shall sketch the proof of Proposition \ref{pr1}.

\section{Existence of minimizers}
The main purpose of this section is to prove Theorem \ref{th1} and Theorem \ref{th2} on the existence and nonexistence of minimizers for $e(a)$. The proofs of Theorem \ref{th1} and Theorem \ref{th2} are based on the following Gagliardo-Nirenberg type inequality in any given open bounded domain.
\begin{prop}\label{pr1}
Suppose $\Omega\subset\mathbb{R}^{N}$ is an open bounded domain, then we have
\begin{equation}\label{2.1}
\frac{a^{*}}{1+\beta^{2}}\int_{\Omega}|x|^{-b}|u(x)|^{2+2\beta^{2}}dx\leq\int_{\Omega}|\nabla u(x)|^{2}dx\Big(\int_{\Omega}|u(x)|^{2}dx\Big)^{\beta^{2}},\ \ u\in H_{0}^{1}(\Omega),
\end{equation}
where the optimal constant $\frac{1+\beta^{2}}{a^{*}}=\frac{1+\beta^{2}}{\|Q\|_{L^{2}(\mathbb{R}^{N})}^{2\beta^{2}}}$ is not attained.
\end{prop}
Since the proof of Proposition \ref{pr1} is long and similar to the proof of \cite[Proposition 2.1]{GL}, we leave it to the Appendix for simplicity. Making full use of Proposition \ref{pr1}, we then complete the proofs of Theorem \ref{th1} and Theorem \ref{th2}.
\vskip 0.1truein
\noindent\textbf{Proof of Theorem \ref{th1}.}
(1). We first proof that $e(a)$ admits at least one minimizer for all $0<a<a^{*}$. Indeed, for any fixed $0<a<a^{*}$ and $u\in\mathcal{M}$, we derive from the nonnegativity of  $V(x)$ and the Gagliardo-Nirenberg inequality (\ref{2.1}) that
\begin{equation}\label{eq2.2}
\begin{split}
E_{a}(u)&=\int_{\Omega}(|\nabla u(x)|^{2}+V(x)|u(x)|^{2})dx-\frac{a}{1+\beta^{2}}\int_{\Omega}|x|^{-b}|u(x)|^{2+2\beta^{2}}dx\\
&\geq (1-\frac{a}{a^{*}})\int_{\Omega}|\nabla u(x)|^{2}dx+\int_{\Omega}V(x)|u(x)|^{2}dx\\
&\geq (1-\frac{a}{a^{*}})\int_{\Omega}|\nabla u(x)|^{2}dx\geq0,
\end{split}
\end{equation}
which implies that $E_{a}(u)$ is bounded uniformly from below.
Choose a minimizing sequence $\{u_{n}\}\subset H_{0}^{1}(\Omega)$ satisfying $\|u_{n}\|_{2}^{2}=1$ and $\lim_{n\to\infty}E_{a}(u_{n})=e(a)$. It then follows from (\ref{eq2.2}) that $\int_{\Omega}|\nabla u_{n}|^{2}dx$ is bounded uniformly for $n$. Since the embedding $H_{0}^{1}(\Omega)\hookrightarrow L^{q}(\Omega)$ $(2\leq q<2^{*})$ is compact, there exists a subsequence, still denoted by $\{u_{n}\}$, such that for some $u\in H_{0}^{1}(\Omega)$,
\begin{equation}\label{eq2.3}
\begin{split}
u_{n}\rightharpoonup u\ \ \text{weakly in}\ \ H_{0}^{1}(\Omega)\ \ \text{and}\ \
u_{n}\rightarrow u\ \ \text{strongly in}\ \ L^{q}(\Omega),\ \ 2\leq q<2^{*}.
\end{split}
\end{equation}
This further implies that $u\in\mathcal{M}$.

As for the convergence of nonlinear term, motivated by \cite{AA}, we claim that
\begin{equation}\label{eq2.4}
\begin{split}
\lim_{n\rightarrow\infty}\int_{\Omega}|x|^{-b}|u_{n}(x)|^{2+2\beta^{2}}dx=
\int_{\Omega}|x|^{-b}|u(x)|^{2+2\beta^{2}}dx.
\end{split}
\end{equation}
Actually, for any $1<\bar{s}<\frac{N}{b}$, applying the H\"{o}lder inequality then yields that
\begin{equation}\label{eq2.5}
\begin{split}
&\Big|\int_{\Omega}|x|^{-b}|u_{n}(x)|^{2+2\beta^{2}}dx-\int_{\Omega}|x|^{-b}|u(x)|^{2+2\beta^{2}}dx\Big|\\
\leq&\int_{\Omega}|x|^{-b}\Big||u_{n}(x)|^{2+2\beta^{2}}-|u(x)|^{2+2\beta^{2}}\Big|dx\\
\leq&\Big(\int_{\Omega}|x|^{-b\bar{s}}dx\Big)^{\frac{1}{\bar{s}}}
\Big(\int_{\Omega}\Big||u_{n}(x)|^{2+2\beta^{2}}-|u(x)|^{2+2\beta^{2}}\Big|^{s}dx\Big)^{\frac{1}{s}}\\
\leq&C\Big(\int_{\Omega}\Big||u_{n}(x)|^{2+2\beta^{2}}-|u(x)|^{2+2\beta^{2}}\Big|^{s}dx\Big)^{\frac{1}{s}},
\end{split}
\end{equation}
where $s, \bar{s}>1$ and $\bar{s}=\frac{s}{s-1}$.
In addition,
\begin{equation}\label{eq2.6}
\begin{split}
&\Big(\int_{\Omega}\Big| |u_{n}|^{2+2\beta^{2}}-|u|^{2+2\beta^{2}}\Big|^{s}dx\Big)^{\frac{1}{s}}\\
\leq& C\Big(\int_{\Omega}|u_{n}-u|^{s}\big(|u_{n}|^{1+2\beta^{2}}+|u|^{1+2\beta^{2}}\big)^{s}dx\Big)
^{\frac{1}{s}}\\
\leq&C\|u_{n}-u\|_{L^{\sigma}(\Omega)}\Big\||u_{n}|^{1+2\beta^{2}}+|u|^{1+2\beta^{2}}\Big\|
_{L^{\frac{\tau}{1+2\beta^{2}}}(\Omega)}\\
\leq&C\|u_{n}-u\|_{L^{\sigma}(\Omega)}\Big(\|u_{n}\|_{L^{\tau}(\Omega)}^{1+2\beta^{2}}+
\|u\|_{L^{\tau}(\Omega)}^{1+2\beta^{2}}\Big),
\end{split}
\end{equation}
where
\begin{equation}\label{eq2.7}
\begin{split}
\frac{1}{\sigma}+\frac{1+2\beta^{2}}{\tau}=\frac{1}{s}<\frac{N-b}{N}.
\end{split}
\end{equation}
On the other hand, for $N=1, 2$, it is obvious that there exist $\sigma, \tau\in[2, \infty)$ such that
$\frac{1}{\sigma}+\frac{1+2\beta^{2}}{\tau}=\frac{1}{s}<\frac{N-b}{N}$; For $N\geq3$, one can check that
$$\frac{N-2}{2N}+\frac{N-2}{2N}(1+2\beta^{2})=\frac{(N-2)(N+2-b)}{N^{2}}<\frac{N-b}{N},$$
which implies that there exist suitable constants $\sigma, \tau\in[2, \frac{2N}{N-2})$ satisfying (\ref{eq2.7}). Combining this with the convergence of $u_{n}$ in (\ref{eq2.3}), we deduce that
$$\|u_{n}\|_{L^{\tau}(\Omega)}^{1+2\beta^{2}}+
\|u\|_{L^{\tau}(\Omega)}^{1+2\beta^{2}}\leq C\ \ \text{and} \ \ \|u_{n}-u\|_{L^{\sigma}(\Omega)}\rightarrow 0\ \ \text{as} \ \ n\rightarrow\infty.$$
From the above facts,
the claim (\ref{eq2.4}) is therefore proved.

Consequently, following the weak lower semicontinuity and (\ref{eq2.4}), we deduce that $$E_{a}(u)\geq e(a)=\lim_{n\to\infty}E_{a}(u_{n})\geq E_{a}(u),$$
which indicates that $u$ is a minimizer of $e(a)$ and Theorem \ref{1.1} (1) is thus proved.
\vskip 0.1truein

(2). We next show that there is no minimizer for $e(a)$ once $a>a^{*}$. Since $\Omega$ is an open bounded domain of $\mathbb{R}^{N}$ and contains $0$, there exists an open ball $B_{2R}(0)\subset\Omega$ centered at an inner point $0$, where $R>0$ is sufficiently small. Choose a nonnegative cutoff function $\varphi(x)\in C_{0}^{\infty}(\mathbb{R}^{N})$ such that $\varphi(x)=1$ for $|x|\leq R$ and $\varphi(x)=0$ for $|x|\geq2R$. Set for all $\tau>0$,
\begin{equation}\label{eq2.8}
\Phi_{\tau}(x)=A_{\tau}\frac{\tau^{\frac{N}{2}}}{\|Q\|_{2}}\varphi(x)Q(\tau x),\ \,x\in\Omega,
\end{equation}
where $A_{\tau}>0$ is chosen such that $\int_{\Omega}|\Phi_{\tau}(x)|^{2}dx=1$. According to the exponential decay of $Q$, we then obtain from (\ref{eq2.8}) that
\begin{equation}\label{eq2.9}
\begin{split}
\frac{1}{|A_{\tau}|^{2}}&=\frac{1}{\|Q\|_{2}^{2}}\int_{\mathbb{R}^{N}}\varphi^{2}(\frac{x}
{\tau})Q^{2}(x)dx
=1+O(\tau^{-\infty})\ \ \text{as}\ \ \tau\rightarrow\infty,
\end{split}
\end{equation}
where 
$O(t^{-\infty}):=m(t)$ denotes the function $m(t)$ satisfying $\lim_{t\rightarrow\infty}|m(t)|t^{s}=0$ for all $s>0$.
Moreover, following (\ref{1.8}) and (\ref{1.9}), some calculations yield that
\begin{equation}\label{eq2.10}
\begin{split}
&\int_{\Omega}|\nabla\Phi_{\tau}(x)|^{2}dx-\frac{a}{1+\beta^{2}}\int_{\Omega}|x|^{-b}
|\Phi_{\tau}(x)|^{2+2\beta^{2}}dx\\
=&\frac{\tau^{2}A_{\tau}^{2}}{\|Q\|_{2}^{2}}\int_{\mathbb{R}^{N}}|\nabla Q(x)|^{2}dx-\frac{a} {1+\beta^{2}}\frac{\tau^{2}A_{\tau}^{2+2\beta^{2}}}{\|Q\|_{2}^{2+2\beta^{2}}}\int_{\mathbb{R}^{N}}|x|^{-b}
|Q(x)|^{2+2\beta^{2}}dx+O(\tau^{-\infty})\\
=&\frac{\tau^{2}A_{\tau}^{2}}{\|Q\|_{2}^{2}}\int_{\mathbb{R}^{N}}|\nabla Q(x)|^{2}dx-\frac{a\tau^{2}A_{\tau}^{2+2\beta^{2}}}{\|Q\|_{2}^{2+2\beta^{2}}}\int_{\mathbb{R}^{N}}|\nabla Q(x)|^{2}dx+O(\tau^{-\infty})\\
=&\frac{\tau^{2}}{\|Q\|_{2}^{2}}(1-\frac{a}{a^{*}})\int_{\mathbb{R}^{N}}|\nabla Q(x)|^{2}dx+O(\tau^{-\infty})\ \ \text{as}\ \ \tau\rightarrow\infty.
\end{split}
\end{equation}
On the other hand, since $x\mapsto V(x)\varphi^{2}(x)$ is bounded in $\Omega$, it follows from \cite[Theorem 1.8]{LE} that
\begin{equation}\label{eq2.11}
\lim_{\tau\rightarrow\infty}\int_{\Omega}V(x)|\Phi_{\tau}(x)|^{2}dx=V(0).
\end{equation}
We then conclude from (\ref{eq2.10}) and (\ref{eq2.11}) that
\begin{equation}\label{eq2.12}
\begin{split}
e(a)&\leq E_{a}(\Phi_{\tau}(x))\\
&=\int_{\Omega}|\nabla\Phi_{\tau}(x)|^{2}dx-\frac{a}{1+\beta^{2}}\int_{\Omega}|x|^{-b}
|\Phi_{\tau}(x)|^{2+2\beta^{2}}dx+\int_{\Omega}V(x)|\Phi_{\tau}(x)|^{2}dx\\
&=\frac{\tau^{2}}{\|Q\|_{2}^{2}}(1-\frac{a}{a^{*}})\int_{\mathbb{R}^{N}}|\nabla Q(x)|^{2}dx+\int_{\Omega}V(x)|\Phi_{\tau}(x)|^{2}dx+O(\tau^{-\infty})\\
&\longrightarrow -\infty\ \ \text{as}\ \ \tau\rightarrow\infty,
\end{split}
\end{equation}
which implies that $e(a)=-\infty$ as soon as $a>a^*$ and the nonexistence of minimizers is therefore proved. This completes the proof of Theorem \ref{th1}.
\qed
\vskip 0.1truein
Next, we shall make full use of the various  Gagliardo-Nirenberg type inequalities to prove Theorem \ref{th2} by contradiction.
\vskip 0.1truein
\noindent\textbf{Proof of Theorem \ref{th2}.}
(1). Applying the Gagliardo-Nirenberg inequality $(\ref{2.1})$, one can check that $e(a^{*})\geq0$ is bounded from below, and thus there exists a minimizing sequence $\{u_{n}\}$ of $e(a^{*})$ such that $e(a^{*})=\lim_{n\rightarrow\infty}E_{a^{*}}(u_{n})$. To prove Theorem \ref{th2} (1), the argument similar to that of proving Theorem \ref{th1} (1),
it is sufficient to prove that $\{u_{n}\}$ is bounded uniformly in $H_{0}^{1}(\Omega)$. On the contrary, suppose that $\{u_{n}\}$ is unbounded in $H_{0}^{1}(\Omega)$, then there exists a subsequence of $\{u_{n}\}$, still denoted by $\{u_{n}\}$, such that $\|u_{n}\|_{H_{0}^{1}(\Omega)}\rightarrow\infty$ as $n\rightarrow\infty$, i.e.,
\begin{equation}\label{eq2.13}
\int_{\Omega}|\nabla u_{n}(x)|^{2}dx\rightarrow\infty\ \ \text{as}\ \ n\rightarrow\infty.
\end{equation}
Define now
\begin{equation}\label{eq2.14}
\varepsilon_{n}^{-2}:=\int_{\Omega}|\nabla u_{n}(x)|^{2}dx,
\end{equation}
then by (\ref{eq2.13}), we have that $\varepsilon_{n}\rightarrow 0$ as $n\rightarrow\infty$. In view of the above facts, we next define
\begin{equation}\label{eq2.15}
w_{n}(x):=\arraycolsep=1.5pt\left\{\begin{array}{lll}
	 \varepsilon_{n}^{\frac{N}{2}}u_{n}(\varepsilon_{n}x)   \quad   &\mbox{if} & \ \, x\in\Omega_{n}:=\{x\in\mathbb{R}^{N}:\varepsilon_{n}x\in\Omega\},\\[4mm]
	\displaystyle 0  \quad   &\mbox{if}& \,\ x\in\mathbb{R}^{N}\setminus\Omega_{n},
\end{array}\right.
\end{equation}
then
\begin{equation}\label{eq2.16}
\int_{\R^N}|\nabla w_{n}(x)|^{2}dx=\int_{\R^N}|w_{n}(x)|^{2}dx=1.
\end{equation}
This implies that $w_n$ is bounded uniformly in $H^1(\R^N)$. According to the embedding theorem, we can extract a subsequence if necessary such that
$$w_{n}\rightharpoonup w_0\,\ \hbox{weakly in} \,\ H^1(\R^N)\ \hbox{and} \  w_{n}\rightarrow w_0\,\ \hbox{strongly in}\,\ L_{loc}^{q}(\R^{N}),\,\ 2\leq q<2^* $$
for some $w_0\in H^1(\R^N)$.

Next, we claim that
\begin{equation}\label{eq2.17}
\lim_{n\to\infty}w_n=\frac{\beta^{\frac{N}{2}}Q(\beta x)}{\|Q\|_2}\quad\hbox{strongly in $H^1(\R^N)$}.
\end{equation}
Indeed, by the definition of $w_n$, we deduce from the Gagliardo-Nirenberg inequality (\ref{1.6}) that
\begin{equation}\label{eq2.18}
\begin{split}
C&\geq e(a^*)=\lim_{n\rightarrow\infty}E_{a^*}(u_{n})
\\&=\lim_{n\rightarrow\infty}\Big\{\int_{\Omega}\Big(|\nabla u_{n}|^{2}-\frac{a^*}{1+\beta^2}|x|^{-b}|u_{n}|^{2+2\beta^2}+V(x)|u_{n}(x)|^{2}\Big)dx\Big\}
\\
&=\lim_{n\rightarrow\infty}\Big\{\frac{1}{\varepsilon^{2}_{n}}\int_{\R^{N}}\Big[|\nabla w_{n}|^{2}-\frac{a^{\ast}}{1+\beta^2}|x|^{-b}|w_{n}|^{2+2\beta^2}\Big]dx+\int_{\Omega_{n}}V(\varepsilon_{n}x)|w_{n}(x)|^{2}dx\Big\}\\
&\geq 0.
\end{split}
\end{equation}
Therefore, we derive from (\ref{eq2.16}), the fact $\eps_n\to 0$ as $n\to\infty$ and above that
\begin{equation}\label{eq2.19}
\lim_{n\to\infty}\frac{a^{\ast}}{1+\beta^2}\int_{\R^{N}}|x|^{-b}|w_{n}(x)|^{2+2\beta^2}dx=\lim_{n\to\infty}\int_{\R^{N}}|\nabla w_{n}(x)|^{2}dx=1.
\end{equation}
We then prove
\begin{equation}\label{eq2.20}
\begin{split}
\lim_{n\rightarrow\infty}\int_{\mathbb{R}^{N}}|x|^{-b}|w_{n}(x)|^{2+2\beta^{2}}dx=
\int_{\mathbb{R}^{N}}|x|^{-b}|w_{0}(x)|^{2+2\beta^{2}}dx.
\end{split}
\end{equation}
Indeed,
\begin{equation}\label{eq2.21}
\begin{split}
&\Big|\int_{\mathbb{R}^{N}}|x|^{-b}|w_{n}(x)|^{2+2\beta^{2}}dx-
\int_{\mathbb{R}^{N}}|x|^{-b}|w_{0}(x)|^{2+2\beta^{2}}dx\Big|\\
\leq&\int_{B_{R(0)}}|x|^{-b}\Big||w_{n}|^{2+2\beta^{2}}-|w_{0}|^{2+2\beta^{2}}\Big|dx+
\int_{\mathbb{R}^{N}\backslash B_{R(0)}}|x|^{-b}\Big||w_{n}|^{2+2\beta^{2}}-|w_{0}|^{2+2\beta^{2}}\Big|dx\\
:=&A_{n}+B_{n}.
\end{split}
\end{equation}
For any $\varepsilon>0$, we can choose $R$ large enough such that
$R> \varepsilon^{-\frac{1}{b}}$, then we derive that
\begin{equation}\label{eq2.22}
\begin{split}
B_{n}&=\int_{\mathbb{R}^{N}\backslash B_{R(0)}}|x|^{-b}\Big||w_{n}|^{2+2\beta^{2}}-|w_{0}|^{2+2\beta^{2}}\Big|dx\leq C\varepsilon.
\end{split}
\end{equation}
On the other hand, the convergence in (\ref{eq2.4}) indicates that
\begin{equation}\label{eq2.23}
\begin{split}
\lim_{n\rightarrow\infty}A_{n}=\lim_{n\rightarrow\infty}
\int_{B_{R(0)}}|x|^{-b}\Big||w_{n}|^{2+2\beta^{2}}-|w_{0}|^{2+2\beta^{2}}\Big|dx=0.
\end{split}
\end{equation}
Combining this with (\ref{eq2.22}), (\ref{eq2.20}) is therefore proved. It then follows from (\ref{eq2.19}) and (\ref{eq2.20}) that
\begin{equation}\label{eq2.24}
\begin{split}
\frac{a^{\ast}}{1+\beta^2}\int_{\mathbb{R}^{N}}|x|^{-b}|w_{0}|^{2+2\beta^{2}}dx=
\lim_{n\rightarrow\infty}\frac{a^{\ast}}{1+\beta^2}\int_{\mathbb{R}^{N}}|x|^{-b}|w_{n}|^{2+2\beta^{2}}dx=1,
\end{split}
\end{equation}
which imples $w_{0}\not\equiv 0$.

Moreover, combining the Gagliardo-Nirenberg inequality (\ref{1.6}) and the fact that $w_{n}\rightharpoonup w_0\,\ \hbox{weakly in} \,\ H^1(\R^N)$, we derive that
\begin{equation}\label{eq2.25}
\begin{split}
0&= \lim_{n\rightarrow\infty}\varepsilon_n^2e(a^*)
\\
&\geq\lim_{n\rightarrow\infty}\int_{\R^{N}}\Big[|\nabla w_{n}(x)|^{2}-\frac{a^{\ast}}{1+\beta^2}|x|^{-b}|w_{n}(x)|^{2+2\beta^2}\Big]dx
\\&\geq\int_{\R^{N}}|\nabla w_{0}(x)|^{2}dx-\frac{a^{\ast}}{1+\beta^2}\int_{\R^{N}}|x|^{-b}|w_{0}(x)|^{2+2\beta^2}dx
\\&\geq\Big[1-\Big(\int_{\R^{N}}|w_{0}(x)|^2dx\Big)^{\beta^2}\Big]\int_{\R^{N}}|\nabla
w_{0}(x)|^{2}dx.
\end{split}
\end{equation}
Since $\|w_0\|_{L^2(\R^N)}\leq 1$, the above inequality implies that
\begin{equation}\label{eq2.26}
\int_{\R^{N}}|w_{0}(x)|^{2}dx=1,\,\ \int_{\R^{N}}|\nabla w_{0}(x)|^{2}dx=\frac{a^{\ast}}{1+\beta^2}\int_{\R^{N}}|x|^{-b}|w_{0}(x)|^{2+2\beta^2}dx=1.
\end{equation}
As a consequence, we conclude that the Gagliardo-Nirenberg inequality (\ref{1.6}) is achieved by $w_0$ and hence $w_{0}=\frac{\beta^{\frac{N}{2}}Q(\beta x)}{\|Q\|_2}$ in view of (\ref{eq2.26}) and (\ref{1.9}).  By the norm preservation, we finally conclude that $w_{n}$ converges to $w_{0}=\frac{\beta^{\frac{N}{2}}Q(\beta x)}{\|Q\|_2}$ strongly in $H^{1}(\R^{N})$ and the claim (\ref{eq2.17}) is thus proved.

Based on (\ref{eq2.17}), applying Fatou's lemma, one can deduce that
\begin{equation}\label{eq2.27}
e(a^{\ast})\geq\lim_{n\rightarrow\infty}\int_{\Omega_{n}}V(\varepsilon_{n} x)|w_{n}(x)|^{2}dx\geq \int_{\R^{N}}V(0)|w_{0}(x)|^{2}dx=V(0),
\end{equation}
which however contradicts the assumption  $e(a^{\ast})<V(0)$. Hence the minimizing sequence $\{u_{n}\}$ of $e(a^*)$ is bounded uniformly in $H_{0}^{1}(\Omega)$ and the proof of Theorem \ref{th2} (1) is accomplished.
\vskip 0.1truein
(2). We next consider the case $a=a^{*}$ and $V(0)=0$. In this case, the argument of proving (\ref{eq2.10}) and (\ref{eq2.11}) yields that
$e(a^{*})\leq V(0)=0$. Additionally, it follows from (\ref{eq2.2}) that $e(a^{*})\geq0$, we then conclude that $e(a^{*})=0$. Suppose now that there exists a minimizer $u$ for $e(a^{*})$,  we would thus have
\begin{equation}\label{eq2.28}
\int_{\Omega}|\nabla u(x)|^{2}dx=\frac{a^{*}}{1+\beta^{2}}\int_{\Omega}|x|^{-b}|u(x)|^{2+2\beta^{2}}dx,
\end{equation}
which contradicts the fact that the optimal constant in Proposition \ref{2.1} is not attained. Therefore there is no minimizer for  $e(a^{*})$ in this case.

Moreover, taking the same test function as in (\ref{eq2.8}), we have $V(0)=\inf_{x\in\Omega}V(x)\leq \lim_{a\nearrow a^{*}}e(a)\leq V(0)$ and hence
$\lim_{a\nearrow a^{*}}e(a)=V(0)=e(a^{*})$. The proof of Theorem \ref{th2} is therefore finished.
\qed
\section{Mass concentration as $a \nearrow a^{*}$}
In this section, we shall analyze the mass concentration behavior of minimizers for $e(a)$ as $a\nearrow a^{*}$ in the case that $V(x)$ satisfies the additional condition $(V_2)$. Let $u_{a}$ be a positive minimizer of $e(a)$, by variational theory,  $u_{a}$ solves the following Euler-Lagrange equation
\begin{equation}\label{3.1}
-\Delta u_{a}(x)+V(x)u_{a}(x)=\mu_a u_{a}(x)+a u_{a}^{1+2\beta^2}(x)|x|^{-b}\ \ \text{in}\ \ \Omega,
\end{equation}
where $\mu_{a}\in\R$ is a suitable Lagrange multiplier satisfying
\begin{equation}\label{3.2}
\mu_a=e(a)-\frac{a\beta^2}{1+\beta^2}\int_{\Omega} |x|^{-b}|u_a|^{2+2\beta^2}dx.
\end{equation}
Motivated by \cite{LL, ZL}, we begin with the following energy estimate.
\begin{lem}\label{lemma3.1}
Under the assumptions of Theorem \ref{th3}, we have the following energy estimate
\begin{equation}\label{3.3}
0\leq e(a)\leq\frac{p+2}{p}\lambda^{2}\Big(\frac{1}{\|Q\|_{2}^{2}}\Big)^{\frac{2}{p+2}}
\Big(\frac{a^*-a}{a^*\beta^{2}}\Big)^{\frac{p}{p+2}} \ \ \mbox{as} \ \ a\nearrow a^*,
\end{equation}
where $p>0$ and $\lambda>0$ are defined by (\ref{1.11}) and (\ref{1.12}) respectively.
\end{lem}
\vskip 0.1truein
\noindent\textbf{Proof.} Since $\Omega$ containing 0 is an open bounded domain of $\mathbb{R}^{N}$, there exists an open ball $B_{2R}(0)\subset\Omega$, where $R>0$ is small enough. Choose the same trial function as in (\ref{eq2.8}), it then follows from (\ref{eq2.9}) and the exponential decay of $Q$ that
\begin{equation*}
\begin{split}
\int_{\Omega}V(x)\Phi_{\tau}^{2}(x)dx
&=\frac{A_{\tau}^2\tau^{N}}{\|Q\|_{2}^2}\int_{B_{2R}(0)}V(x)\varphi^{2}(x)Q^2(\tau x)dx\\
&=\frac{A_{\tau}^2}{\|Q\|_{2}^2}\int_{B_{2\tau R}(0)}V(\frac{x}{\tau})\varphi^{2}(\frac{x}{\tau})Q^2(x)dx\\
&\leq\frac{A_{\tau}^2}{\|Q\|_{2}^2}\int_{B_{\sqrt{\tau R}}(0)}V(\frac{x}{\tau})Q^2( x)dx+Ce^{-\sqrt{\tau R}}\\
&=\frac{2}{p\|Q\|_{2}^{2}}\lambda^{p+2}\tau^{-p}+o(\tau^{-p})\ \ \ \text{as}\ \ \tau\rightarrow\infty.
\end{split}
\end{equation*}
Therefore, we derive from (\ref{1.9}), (\ref{eq2.10}) and above that
\begin{equation*}
\begin{split}
e(a)\leq E_{a}(\Phi_{\tau})&=\int_{\Omega}|\nabla\Phi_{\tau}(x)|^{2}dx-\frac{a}{1+\beta^{2}}\int_{\Omega}|x|^{-b}
|\Phi_{\tau}(x)|^{2+2\beta^{2}}dx+\int_{\Omega}V(x)|\Phi_{\tau}(x)|^{2}dx\\
&=\frac{\tau^{2}}{\|Q\|_{2}^{2}}(1-\frac{a}{a^{*}})\int_{\mathbb{R}^{N}}|\nabla Q(x)|^{2}dx+\int_{\Omega}V(x)|\Phi_{\tau}(x)|^{2}dx+O(\tau^{-\infty})\\
&\leq\frac{(a^*-a)\tau^{2}}{a^*\|Q\|_{2}^2}\int_{\R^{N}}|\nabla Q(x)|^{2}dx+\frac{2}{p\|Q\|_{2}^2}\lambda^{p+2}\tau^{-p}+o(\tau^{-p})\\
&=\frac{(a^*-a)\tau^{2}}{a^*\beta^{2}}+\frac{2}{p\|Q\|_{2}^2}\lambda^{p+2}\tau^{-p}+o(\tau^{-p})
\ \ \ \text{as}\ \ \tau\rightarrow\infty.
\end{split}
\end{equation*}
By choosing $\tau=\lambda\big[\frac{a^*\beta^{2}}{(a^*-a)\|Q\|_{2}^2}\big]^{\frac{1}{p+2}}$, one can check that
$$ e(a)\leq\frac{p+2}{p}\lambda^{2}\Big(\frac{1}{\|Q\|_{2}^{2}}\Big)^{\frac{2}{p+2}}
\Big(\frac{a^*-a}{a^*\beta^{2}}\Big)^{\frac{p}{p+2}} \ \ \mbox{as} \ \ a\nearrow a^*. $$
Combining this with (\ref{eq2.2}), we complete the proof of Lemma \ref{lemma3.1}.
\qed
\vskip 0.1truein
We then establish the following crucial lemma, which is a weak version of Theorem \ref{th3}.
\vskip 0.1truein
\begin{lem}\label{lemma3.2}
Under the assumptions of Theorem \ref{th3}, and suppose $u_{a}$ be a positive minimizer of $e(a)$. Define
\begin{equation}\label{3.4}
\varepsilon_{a}:=\Big(\int_{\Omega}|\nabla u_{a}|^{2}dx\Big)^{-\frac{1}{2}},
\end{equation}
and
\begin{equation}\label{3.5}
w_{a}(x):=\varepsilon^{\frac{N}{2}}_{a}u_{a}(\varepsilon_{a}x),\ \ x\in \Omega_{a},
\end{equation}
where $\Omega_{a}:=\{x\in\mathbb{R}^{N}:\varepsilon_{a}x\in\Omega\}$.
We then have the following
\begin{enumerate}
\item $\varepsilon_{a}\rightarrow 0$ as $a\nearrow a^{\ast}$ and $\lim_{a\nearrow a^{\ast}}\Omega_{a}=\mathbb{R}^{N}$.
\item $w_a(x)\rightarrow\frac{\beta^{\frac{N}{2}}Q(\beta x)}{\|Q\|_2}$ strongly in $ H^1(\R^N)\cap L^{\infty}(\R^{N})$ as $a\nearrow a^{\ast}$ in the sense that $w_{a}(x)\equiv0$ for all $x\in\mathbb{R}^{N}\setminus\Omega_{a}$.
\item There exist positive constants $R_0>0$ and $C>0$ independent of $a$ such that
    \begin{equation}\label{f1}
|w_{a}(x)|, \ \ |\nabla w_a(x)|\leq Ce^{-\frac{\beta|x|}{2}}\ \ \mbox{in}\ \ \Omega_{a}\backslash B_{R_{0}}(0)\ \ \mbox{as}\ \ a\nearrow a^*.
\end{equation}
\end{enumerate}
\end{lem}
\noindent\textbf{Proof.}  1. By contradiction, assume that there exists a sequence $\{a_{k}\}$ satisfying $a_{k}\nearrow a^{\ast}$ as $k\rightarrow\infty$ such that the sequence $\{u_{a_{k}}\}$ is bounded uniformly in $H_{0}^1(\Omega)$. By applying the compact embedding
$H_{0}^{1}(\Omega)\hookrightarrow L^{q}(\Omega)$ for $2\leq q<2^*$, we can choose a subsequence of $\{u_{a_{k}}\}$, still denoted by $\{u_{a_{k}}\}$, and $u_0\in H_{0}^1(\Omega)$ such that
$${u_{a_{k}}}\rightharpoonup u_{0}\,\ \hbox{weakly in}\,\ \ H_{0}^1(\Omega)\,\ \hbox{and}\,\   {u_{a_{k}}}\rightarrow u_{0}\,\ \hbox{strongly in}\,\  L^{q}(\Omega),\,\ 2\leq q <2^*.$$
Combining the above convergence with (\ref{eq2.4}) gives that
\begin{equation*}
0=e(a^{\ast})\leq E_{a^{\ast}}(u_{0})\leq \lim\limits_{k\rightarrow\infty}E_{a_{k}}(u_{a_{k}})=\lim\limits_{k\rightarrow\infty}e(a_{k})=0,
\end{equation*}
which indicates that $u_{0}$ is a minimizer of $e(a^*)$. This however contradicts Theorem \ref{th2} (2), which further implies that $\varepsilon_{a}\rightarrow 0$ as  $a\nearrow a^{\ast}$. By the definition of $\Omega_{a}$, we deduce that $\lim_{a\nearrow a^{\ast}}\Omega_{a}=\mathbb{R}^{N}$ and the proof of Lemma \ref{lemma3.2} (1) is complete.
\vskip 0.1truein

2. Note from Theorem \ref{th2} (2) that $0=e(a^*)=\lim_{a\nearrow a^*}e(a)$, thus $\{u_a\}$ is also a minimizing sequence of $e(a^*)$.
In a similar way to (\ref{eq2.15}), we define the zero continuation of $w_{a}(x)$ as follows
\begin{equation}\label{3.6}
\tilde{w}_{a}(x):=\arraycolsep=1.5pt\left\{\begin{array}{lll}
	 w_{a}(x)   \quad   &\mbox{if} & \ \, x\in\Omega_{a},\\[4mm]
	\displaystyle 0  \quad   &\mbox{if}& \,\ x\in\mathbb{R}^{N}\setminus\Omega_{a}.
\end{array}\right.
\end{equation}
Since $\lim_{a\nearrow a^{*}}\Omega_{a}=\mathbb{R}^{N}$, combining this with the exponential decay of $Q(x)$, in order to prove Lemma \ref{lemma3.2} (2), we just need to prove that
$$\tilde{w}_a(x)\rightarrow\frac{\beta^{\frac{N}{2}}Q(\beta x)}{\|Q\|_2}\ \ \text{strongly in} \ \ H^1(\R^N)\cap L^{\infty}(\R^{N})\ \ \text{as} \ \ a\nearrow a^{\ast}.$$

Firstly, we can derive from(\ref{3.4}), (\ref{3.5}) and (\ref{3.6}) that
\begin{equation}\label{7}
\int_{\R^{N}} |\nabla \tilde{w}_a|^2dx=\int_{\R^{N}} |\tilde{w}_a|^2dx=1\,\ \hbox{for any $a\nearrow a^*$.}
\end{equation}
Hence $\tilde{w}_a$ is bounded uniformly for $a$ in $H^1(\R^N)$. Passing to a subsequence if necessary, we obtain that
$$\tilde{w}_{a}\rightharpoonup \tilde{w}_0\,\ \hbox{weakly in} \,\ H^1(\R^N)\ \hbox{and} \  \tilde{w}_{a}\rightarrow \tilde{w}_0\,\ \hbox{strongly in}\,\ L_{loc}^{q}(\R^{N}),\,\ 2\leq q<2^* $$
holds for some $\tilde{w}_0\in H^1(\R^N)$.
Note that the proof of the claim (\ref{eq2.17}) does not rely on the assumption $0\leq e(a^{*})< V(0)$. As a result, the similar argument of proving (\ref{eq2.17}) yields that
\begin{equation}\label{3.7}
\tilde{w}_{a}(x)\rightarrow \frac{\beta^{\frac{N}{2}}Q(\beta x)}{\|Q\|_2}\,\ \hbox{strongly in}\,\ H^1(\R^N)\,\ \hbox{as}\,\ a\nearrow a^*.
\end{equation}

Next, we claim that
\begin{equation}\label{3.8}
\tilde{w}_a(x)\rightarrow \frac{\beta^{\frac{N}{2}}Q(\beta x)}{\|Q\|_2}\,\ \hbox{strongly in}\,\ L^\infty(\R^N)\,\ \hbox{as}\,\  a\nearrow a^*.
\end{equation}
For $N=1$, the above claim is trivial due to the fact that the embedding $H^1(\R)\hookrightarrow L^\infty(\R)$ is continuous. We then consider the case $N\geq2$.
We first prove that
\begin{equation}\label{3.9}
\tilde{w}_{a}(x)\rightarrow 0 \ \ \mbox{as}\ \ |x|\rightarrow\infty \ \ \text{uniformly for}\ \ a\nearrow a^*.
\end{equation}
In fact, we derive from (\ref{7}) that for sufficiently large $R>0$,
\begin{equation*}
\int_{\mathbb{R}^{N}\setminus B_{R}(0)}|\tilde{w}_{a}(x)|^{2}dx\rightarrow0\ \ \text{as} \ \ a\nearrow a^*.
\end{equation*}
Combining this with the definition of $\tilde{w}_{a}(x)$, it then yields from Lemma \ref{lemma3.2} (1) that
\begin{equation}\label{3.10}
\int_{\Omega_{a}\setminus B_{R}(0)}|w_{a}(x)|^{2}dx\rightarrow0 \ \ \text{as}\ \ a\nearrow a^*.
\end{equation}
On the other hand, recall from (\ref{3.1}) and (\ref{3.5}) that $w_{a}$ satisfies the following equation
\begin{equation}\label{3.11}
-\Delta w_a+\eps_a^2V(\eps_ax)w_a=\eps_a^2\mu_aw_a+a w_a^{1+2\beta^2} |x|^{-b}\ \ \hbox{in}\,\ \Omega_{a},
\end{equation}
where $\mu_a\in\R$ is a suitable lagrange multiplier satisfying (\ref{3.2}).
In view of (\ref{3.2}) and the proof of (\ref{eq2.19}), we can obtain that
\begin{equation}\label{3.13}
\begin{split}
\eps_a^2\mu_a&=\eps_a^2e(a)-\eps_a^2\frac{a\beta^2}{1+\beta^2}\int_{\Omega} |x|^{-b}|u_a(x)|^{2+2\beta^2}dx\\
&=\eps_a^2e(a)-\frac{a\beta^2}{1+\beta^2}\int_{\Omega_{a}} |x|^{-b}|w_a(x)|^{2+2\beta^2}dx\\
&=\eps_a^2e(a)-\frac{a\beta^2}{1+\beta^2}\int_{\R^{N}} |x|^{-b}|\tilde{w}_a(x)|^{2+2\beta^2}dx\to -\beta^2\,\ \hbox{as}\,\ a\nearrow a^*.
\end{split}
\end{equation}
It then follows from (\ref{3.11}) and (\ref{3.13}) that
\begin{equation}\label{3.14}
-\Delta w_{a}(x)-c(x)w_{a}(x)\leq 0\ \ \text{in} \ \ \Omega_{a} \ \ \text{as} \ \ a\nearrow a^*,
\end{equation}
where $c(x)=aw_{a}^{2\beta^{2}}|x|^{-b}$. Similar to the proof of Lemma B.2 in \cite{LZ}, one can check that $c(x)\in L^{m}(\Omega_{a})$ where $m\in\Big(\frac{N}{2},\frac{N^2}{2N+2b-4}\Big)$ if $N\geq 3$ and $m\in(1,\frac{2}{b})$ if $N=2$. 
Applying De Giorgi-Nash-Moser theory  \cite[Theorem 4.1]{HL} to (\ref{3.14}), we conclude that there exists $C>0$ independent of $a$ such that
\begin{equation}\label{3.15}
\max\limits_{B_{l}(\rho)} w_{a}(x)\leq C\Big(\int_{{B_{2l}(\rho)}}|w_{a}(x)|^{2} dx\Big)^\frac{1}{{2}}\ \ \mbox{as}\ \ a\nearrow a^*,
\end{equation}
where $\rho\in\Omega_{a}$ is arbitrary and $l>0$ is small enough satisfying $B_{2l}(\rho)\subset\Omega_{a}$. In view of (\ref{3.10}), (\ref{3.15}) and the definition of $\tilde{w}_{a}(x)$, (\ref{3.9}) is therefore proved .

Following the above argument and the exponential decay of $Q(x)$, to complete the proof of (\ref{3.8}), it is enough to show that
\begin{equation}\label{3.16}
\tilde{w}_{a}(x)\rightarrow \frac{\beta^{\frac{N}{2}}Q(\beta x)}{\|Q\|_2} \ \ \mbox{strongly in}\ \ L_{loc}^{\infty}(\R^N)\ \ \mbox{as}\ \  a\nearrow a^*.
\end{equation}
Note that $\lim_{a\nearrow a^{*}}\Omega_{a}=\mathbb{R}^{N}$, for any $R>0$, it holds $B_{R+1}(0)\subset\Omega_{a}$ as $a\nearrow a^{*}$.
 Recall that $w_{a}$ satisfies (\ref{3.11}), we denote
\begin{equation*}
G_a(x):=\varepsilon_a^2\mu_aw_a(x)-\varepsilon_a^2
V(\varepsilon_ax)w_a(x)+
aw^{1+2\beta^2}_{a}|x|^{-b},
\end{equation*}
then
\begin{equation}\label{3.17}
-\Delta w_a(x)=G_a(x) \ \ \mbox{in}\ \ \Omega_{a}.
\end{equation}
Since $\int_{\Omega_{a}}|w_a(x)|^{2}dx=1$, it follows from (\ref{3.15}) that
\begin{equation}\label{3.18}
w_a(x)\ \ \text{is bounded uniformly  in}\ \  L^{\infty}(\Omega_{a}).
\end{equation}
Therefore, for any $r\in(\frac{N}{2},\frac{N}{b})$,
$w_{a}^{1+2\beta^{2}}|x|^{-b}$ is bounded uniformly in $L^{r}(B_{R+1}(0))$,
which implies that $G_{a}(x)$ is bounded uniformly in $L^{r}(B_{R+1}(0))$. It then follows from
\cite[Theorem 9.11]{GT} that
\begin{equation}\label{3.19}
\|w_a(x)\|_{W^{2,r}(B_{R}(0))}\leq C\Big(\|w_a(x)\|_{L^r(B_{R+1}(0))}+\|G_a(x)\|_{L^r(B_{R+1}(0))}\Big),
\end{equation}
where $C>0$ is independent of $a$ and $R$. Therefore, we can obtain that $w_{a}(x)$ is also bounded uniformly in $W^{2,r}(B_{R}(0))$.
Since $r>\frac{N}{2}$ and the embedding $W^{2,r}(B_{R}(0))\hookrightarrow L^{\infty}(B_{R}(0))$ is compact, see \cite[Theorem 7.26]{GT}, we deduce that there exists a subsequence $\{w_{a_{k}}\}$ of $\{w_a\}$ such that
\begin{equation}\label{3.20}
\lim\limits_{a_k\nearrow a^*}\tilde{w}_{a_k}(x)=\lim\limits_{a_k\nearrow a^*}w_{a_k}(x)=\hat{w}_0(x)\ \ \text{strongly in}\ \ L^\infty(B_{R}(0)).
\end{equation}
In view of (\ref{3.7}) and the fact that $R > 0$ is arbitrary, we deduce that
\begin{equation}\label{3.21}
\lim_{a_{k}\nearrow a^*}\tilde{w}_{a_{k}}(x)=\frac{\beta^{\frac{N}{2}}Q(\beta x)}{\|Q\|_2} \ \ \mbox{strongly in}\ \ L_{loc}^{\infty}(\R^N).
\end{equation}
Moreover, because the above convergence is independent of what subsequence we choose, we conclude that $\tilde{w}_{a}(x)\rightarrow\frac{\beta^{\frac{N}{2}}Q(\beta x)}{\|Q\|_2}$ in $ L_{loc}^{\infty}(\R^N)$ as $a \nearrow a^*$ and \eqref{3.16} is proved. Combining this with \eqref{3.9} that \eqref{3.8} holds true.
As a result, by  the definition of $\tilde{w}_{a}(x)$, we derive from (\ref{3.7}) and (\ref{3.8}) that Lemma \ref{lemma3.2} (2) is proved.
\vskip 0.1truein
3. Following \eqref{3.9}, \eqref{3.11} and \eqref{3.13}, we obtain that there exists $R_{0}>0$ large enough such that
\begin{equation}\label{3.22}
-\Delta w_{a}(x)+\frac{\beta^2}{4}w_{a}(x)\leq 0\ \ \mbox{in}\ \ \Omega_{a}\backslash B_{R_{0}}(0) \ \ \mbox{as}\ \ a\nearrow a^*.
\end{equation}
By the comparison principle, comparing $w_{a}(x)$ with $Ce^{-\frac{\beta|x|}{2}}$, we get that
\begin{equation}\label{3.23}
w_{a}(x)\leq Ce^{-\frac{\beta|x|}{2}}\ \ \mbox{in}\ \ \Omega_{a}\backslash B_{R_{0}}(0)\ \ \mbox{as}\ \ a\nearrow a^*,
\end{equation}
where $C>0$ is independent of $a$. Besides, applying the local elliptic estimate \cite[(3.15)]{GT} to \eqref{3.17} yields that
\begin{equation}\label{3.24}
|\nabla w_{a}(x)|\leq Ce^{-\frac{\beta|x|}{4}}\ \ \mbox{in}\ \ \Omega_{a}\backslash B_{R_{0}}(0) \ \ \mbox{as}\ \ a\nearrow a^*.
\end{equation}
Therefore, the exponential decay (\ref{f1}) is then proved and we complete the proof of Lemma \ref{lemma3.2}.
 \qed
~\\

\noindent\textbf{Proof of Theorem 1.3.} By the definition of $w_{a}(x)$, some direct calculations yield that
\begin{equation}\label{3.25}
\begin{split}
e(a)&=E_{a}(u_{a})
\\&=\int_{\Omega}|\nabla u_{a}(x)|^{2}dx-\frac{a}{1+\beta^2}\int_{\Omega}|x|^{-b}u_{a}^{2+2\beta^2}(x)dx+
\int_{\Omega}V(x)|u_{a}(x)|^{2}dx\\
&=\frac{1}{\varepsilon^{2}_{a}}\Big(\int_{\Omega_{a}}|\nabla w_{a}(x)|^{2}dx-\frac{a^{\ast}}{1+\beta^2}\int_{\Omega_{a}}|x|^{-b}w_{a}^{2+2\beta^2}(x)dx\Big)
\\&+\frac{a^{\ast}-a}{(1+\beta^2)\varepsilon^{2}_{a}}\int_{\Omega_{a}}|x|^{-b}w_{a}^{2+2\beta^2}(x)dx
+\int_{\Omega_{a}}V(\varepsilon_{a}x)|w_{a}(x)|^{2}dx\\
&\geq\frac{a^{\ast}-a}{(1+\beta^2)\varepsilon^{2}_{a}}\int_{\Omega_{a}}|x|^{-b}w_{a}^{2+2\beta^2}(x)dx
+\int_{\Omega_{a}}V(\varepsilon_{a}x)|w_{a}(x)|^{2}dx.
\end{split}
\end{equation}
Based on Lemma \ref{lemma3.2} (2) and (\ref{3.23}), one can check that
\begin{equation}\label{3.26}
\begin{split}
&\int_{\Omega_{a}}V(\varepsilon_{a}x)|w_{a}(x)|^{2}dx\\
=&\int_{B_{\varepsilon_{a}^{-\frac{1}{2}}}(0)}V(\varepsilon_{a}x)|w_{a}(x)|^{2}dx+
\int_{\Omega_{a}\setminus B_{\varepsilon_{a}^{-\frac{1}{2}}}(0)}V(\varepsilon_{a}x)|w_{a}(x)|^{2}dx
\\=&\int_{B_{\varepsilon_{a}^{-\frac{1}{2}}}(0)}h(\varepsilon_{a}x)|\varepsilon_{a}x|^{p}
\prod_{i=1}^{n-1}|\varepsilon_{a}x-x_{i}|^{p_{i}}|w_{a}(x)|^{2}dx+o(\varepsilon_{a}^{p})\\
=&\Big(\frac{\varepsilon_{a}}{\beta}\Big)^{p}\frac{1}{\|Q\|_{2}^{2}}
\int_{\mathbb{R}^{N}}|x|^p Q^{2}(x)dx\lim_{x\rightarrow0}\Big(h(x)\prod_{i=1}^{n-1}|x
-x_{i}|^{p_{i}}\Big)+o(\varepsilon_{a}^{p})\\
=&\Big(\frac{\varepsilon_{a}}{\beta}\Big)^{p}\frac{2\lambda^{p+2}}{p\|Q\|_{2}^{2}}
+o(\varepsilon_{a}^{p})\ \ \text{ as} \ \ a \nearrow a^*.
\end{split}\end{equation}
Using the similar argument of the proof for \eqref{eq2.19}, we can deduce that
\begin{equation}\label{3.27}
\begin{split}
&\frac{a^{\ast}-a}{(1+\beta^2)\varepsilon^{2}_{a}}\int_{\Omega_{a}}|x|^{-b}w_{a}^{2+2\beta^2}(x)dx
=[1+o(1)]\frac{(a^*-a)}{a^*\varepsilon^{2}_{a}}\ \ \mbox{as}\ \ a \nearrow a^*.
\end{split}
\end{equation}

In view of the above facts, we conclude that
\begin{equation}\label{3.28}
\begin{split}
e(a)&\geq\frac{a^*-a}{a^*\varepsilon^{2}_{a}}+\Big(\frac{\varepsilon_{a}}{\beta}\Big)^{p}\frac{2
\lambda^{p+2}}{p\|Q\|_{2}^{2}}+o\Big(\varepsilon_{a}^{p}+\frac{a^*-a}{\varepsilon^{2}_{a}}\Big)\\
&\geq\frac{p+2}{p}\lambda^{2}\Big(\frac{1}{\|Q\|_{2}^{2}}\Big)^{\frac{2}{p+2}}
\Big(\frac{a^*-a}{a^*\beta^2}\Big)^{\frac{p}{p+2}}+o\Big((a^*-a)^{\frac{p}{p+2}}\Big)\ \ \mbox{as}\ \ a \nearrow a^*.
\end{split}
\end{equation}
Combining this with the upper energy estimate (\ref{3.3}), we obtain that
\begin{equation}\label{3.29}
\begin{split}
\lim_{a\nearrow a^*}\frac{e(a)}{(a^*-a)^{\frac{p}{p+2}}}=\frac{p+2}{p}\lambda^{2}
\Big(\frac{1}{\|Q\|_{2}^{2}}\Big)^{\frac{2}{p+2}}\Big(\frac{1}{a^*\beta^2}\Big)^{\frac{p}{p+2}},
\end{split}
\end{equation}
and the second equality of (\ref{3.28}) holds if and only if
\begin{equation}\label{3.30}
\begin{split}
\lim_{a\nearrow a^*}\frac{\varepsilon_{a}}{(a^*-a)^{{\frac{1}{p+2}}}}=\frac{1}{\lambda}\Big(\frac{\|Q\|_{2}
^{2}\beta^{p}}{a^{*}}\Big)^{\frac{1}{p+2}}.
\end{split}
\end{equation}
Applying Lemma \ref{lemma3.2}, we then conclude that
\begin{equation}\label{3.31}
\lim_{a\nearrow a^{*}}
\Big(\frac{(a^*-a)\|Q\|_{2}^{2}}{a^*\beta^{2}\lambda^{p+2}}\Big)^{\frac{N}{2(p+2)}}
u_{a}\Big(\Big(\frac{(a^*-a)\|Q\|_{2}^{2}}{a^*\beta^{2}}\Big)^{\frac{1}{p+2}}
\frac{x}{\lambda}
\Big)=\frac{Q(x)}{\|Q\|_{2}}
\end{equation}
strongly in $H^1(\R^N)\cap L^{\infty}(\R^{N})$ and the proof of Theorem \ref{th3} is thus complete.\qed

\section{Local uniqueness of positive minimizers}
Following the $L^{\infty}$-uniform convergence and exponential decay of previous section, this section is concerned with the analysis of the local uniqueness of positive minimizers for $e(a)$ as $a\nearrow a^*$. 
Here and in the sequel, we always denote $u_{k}$ to be a positive minimizer of $e(a_{k})$ for convenience, where $a_{k}\nearrow a^*$ as $k\rightarrow\infty$. By the variational theory, there exists a lagrange multiplier $\mu_{k}\in\mathbb{R}$ satisfying
\begin{equation}\label{4.2}
\mu_{k}=e(a_{k})-\frac{a_{k}\beta^2}{1+\beta^2}\int_{\Omega}|x|^{-b}|u_{k}|^{2+2\beta^{2}}dx,
\end{equation}
such that $u_{k}$ solves the following Euler-Lagrange equation
\begin{equation}\label{4.3}
-\Delta u_{k}(x)+V(x)u_{k}(x)=\mu_{k}u_{k}(x)+a_{k}u_{k}^{1+2\beta^{2}}(x)|x|^{-b}\ \ \text{in} \ \ \Omega.
\end{equation}

Under the assumptions of Theorem \ref{th4}, we also define
\begin{equation}\label{4.4}
\alpha_{k}:=\frac{1}{\lambda}\Big(\frac{(a^{*}-a_{k})\|Q\|_{2}^{2}\beta^p}{a^{*}}\Big)^
{\frac{1}{p+2}}>0.
\end{equation}
It follows from (\ref{3.30}) that
\begin{equation}\label{4.5}
\lim_{k\rightarrow\infty}\frac{\eps_{a_{k}}}{\alp_k}=1,
\end{equation}
which combined with \eqref{3.13} yields that
\begin{equation}\label{4.6}
\mu_{k}\alp_{k}^{2}\rightarrow -\beta^{2}\,\ \mbox{as} \,\ k\rightarrow \infty.
\end{equation}
Based on the above argument, we now prove the local uniqueness of positive minimizers for $e(a)$ as $a\nearrow a^*$.
\vskip 0.1truein
\noindent\textbf{Proof of Theorem 1.4.} Suppose that there are two different positive minimizers $u_{1,k}$ and $u_{2,k}$ of $e(a_{k})$, where $a_{k}\nearrow a^*$ as $k\rightarrow\infty$. It then follows from \eqref{4.3} that $u_{i,k}$ solves the following equation
\begin{equation}\label{4.10}
-\Delta u_{i,k}(x)+V(x)u_{i,k}(x)=\mu_{i,k}u_{i,k}(x)+a_{k}u_{i,k}^{1+2\beta^{2}}(x)|x|^{-b}\ \ \text{in} \ \ \Omega,\  \ i=1, 2,
\end{equation}
where $\mu_{i,k}\in\mathbb{R}$ is the suitable lagrange multiplier. Define
\begin{equation}\label{4.11}
\bar{u}_{i,k}(x):=\frac{\alpha_{k}^{\frac{N}{2}}\|Q\|_{2}}{\beta^{\frac{N}{2}}}\
u_{i,k}(x)\ \ \text{and} \ \
\tilde{u}_{i,k}(x):=\bar{u}_{i,k} \big(\frac{\alp _k}{\beta}x\big),\ \ i=1,2.
\end{equation}
Note from (\ref{4.3}) that $\bar{u}_{i,k}$ satisfies
\begin{equation}\label{4.12}
-\alpha_{k}^{2}\Delta\bar{u}_{i,k}(x)+\alpha_{k}^{2}V(x)\bar{u}_{i,k}(x)=\mu_{i,k}\alpha_{k}^{2}
\bar{u}_{i,k}(x)+\frac{a_{k}\alpha_{k}^{b}\beta^{2-b}}{a^{*}}|x|^{-b}\bar{u}_{i,k}^
{1+2\beta^{2}}(x)\ \ \text{in}\ \ \Omega,
\end{equation}
and $\tilde{u}_{i,k}(x)$ satisfies
\begin{equation}\label{4.8}
-\Delta\tilde{u}_{i,k}(x)+\frac{\alpha_{k}^{2}}{\beta^2}V(\frac{\alpha_{k}}{\beta}x)\tilde{u}_{i,k}(x)
=\frac{\mu_{k}\alpha_{k}^{2}}{\beta^2}\tilde{u}_{i,k}(x)+\frac{a_{k}}{a^{*}}|x|^{-b}\tilde{u}_{i,k}^
{1+2\beta^{2}}(x)\ \ \text{in}\ \ \Omega_{k},
\end{equation}
where $\Omega_{k}:=\{x\in\mathbb{R}^{N}:\frac{\alpha_{k}}{\beta}x\in\Omega\}\rightarrow\mathbb{R}^{N}$ as $k\rightarrow\infty$ in view of Lemma \ref{lemma3.2} (1) and (\ref{4.5}). From now on, when necessary we shall extend $\tilde{u}_{i,k}(x)$ to $\mathbb{R}^{N}$ by setting $\tilde{u}_{i,k}(x)\equiv0$ on $\mathbb{R}^{N}\setminus\Omega_{k}$. It then follows from Section 3 that $\tilde{u}_{i,k}(x)$ is bounded uniformly in $L^{\infty}(\mathbb{R}^{N})$ and $\tilde{u}_{i,k}(x)\rightarrow Q(x)$ uniformly in $\mathbb{R}^{N}$ as $k\rightarrow\infty$.
Moreover, by the exponential decay (\ref{f1}), one can check that there exist $C>0$ and $R_{1}>0$ such that
\begin{equation}\label{4.9}
|\tilde{u}_{i,k}(x)|,\ \ |\nabla\tilde{u}_{i,k}(x)|\leq C e^{-\frac{|x|}{4}}\ \ \mbox{in}\ \ \Omega_{k}\backslash B_{R_{1}}(0)\ \ \mbox{as}\ \ k\rightarrow\infty,
\end{equation}
where $C>0$ is independent of $k$.

Since $u_{1,k}\not\equiv u_{2,k}$, we consider
\begin{equation}\label{4.13}
\bar \xi_k(x)=\frac{ u_{2,k}(x)-  u_{1,k}(x)}{\|  u_{2,k}-  u_{1,k}\|_{L^\infty(\Omega)}}=\frac{\bar u_{2,k}(x)-\bar u_{1,k}(x)}{\|\bar u_{2,k}-\bar u_{1,k}\|_{L^\infty(\Omega)}}.
\end{equation}
Stimulated by \cite{PLC}, we first claim that for any $x_0\in\Omega$, there exists a small constant $\delta >0$  such that
\begin{equation}\label{4.14}
\int_{\partial B_\delta (x_0)} \Big[ \alp ^2_k |\nabla \bar \xi_k|^2+ \frac{\beta ^2}{2} |\bar\xi_k|^2+ \alp ^2_k  V(x)|\bar\xi_k|^2\Big]dS=O( \alp ^N_k)\quad\text{as}\ \ k\rightarrow\infty.
\end{equation}
In fact, denote
\begin{equation}\label{4.15}
\arraycolsep=1.5pt\begin{array}{lll}
\bar{D}^{s-1}_{k}(x)&&:=\displaystyle\frac{\bar{u}_{2,k}^{s}(x)-\bar{u}_{1,k}^{s}(x)}{s(\bar{u}_
{2,k}-\bar{u}_{1,k})}\\[4mm]
&&\displaystyle=\frac{\int^{1}_{0}\frac{d}{dt}[t\bar{u}_{2,k}+(1-t)\bar{u}_{1,k}]^{s}dt}{s(\bar{u}_
{2,k}-\bar{u}_{1,k})}=\int^{1}_{0}\big[t\bar{u}_{2,k}+(1-t)\bar{u}_{1,k}\big]^{s-1}dt,
\end{array}
\end{equation}
then we obtain from (\ref{4.12}) and (\ref{4.13}) that $\bar \xi_k$ satisfies
\begin{equation}\label{4.16}
-\alp ^2_k\Delta \bar \xi_k +\bar C_{k}(x)\bar \xi_k =\bar g_k(x)\quad \text{in}\ \ \Omega,
\end{equation}
where the coefficients
\begin{equation}\label{4.17}
\bar C_k(x):=-\mu _{1,k}\alp ^2_k-\frac{a_{k}\alp^{b}_{k}\beta^{2-b}(1+2\beta^2)}{a^*}|x|^{-b}\bar{D}^{2\beta^2}_{k}+ \alp ^2_k  V(x),
\end{equation}
and
\begin{equation}\label{4.18}
\arraycolsep=1.5pt\begin{array}{lll}
\bar g_k(x)&:=\displaystyle \frac{\bar u_{2,k}(\mu _{2,k}-\mu _{1,k})\alp ^2_k}{\|\bar u_{2,k}-\bar u_{1,k}\|_{L^\infty(\Omega)}}\\[4mm]
&=-\displaystyle\frac{a_{k}\bar u_{2,k}\alpha_{k}^{2}\beta^2}{(1+\beta^2)\|\bar u_{2,k}-\bar u_{1,k}\|_{L^\infty(\Omega)}}\int_{\Omega}|x|^{-b}(u_{2,k}^{2+2\beta^2}-u_{1,k}^{2+2\beta^2})dx\\[4mm]
&=-\displaystyle\frac{a_{k}\bar u_{2,k}\alpha_{k}^{2}\beta^2}{(1+\beta^2)\|\bar u_{2,k}-\bar u_{1,k}\|_{L^\infty(\Omega)}}\frac{\beta^{N(1+\beta^{2})}}{\alpha_{k}^{N(1+\beta^2)}\|Q\|^{2(1+\beta^2)}}
\int_{\Omega}|x|^{-b}(\bar u_{2,k}^{2+2\beta^2}-\bar u_{1,k}^{2+2\beta^2})dx\\[4mm]
&=-\displaystyle  \frac{2a_{k}\bar u_{2,k}\beta^{N+4-b}}{\alp^{N-b}_k\|Q\|^{2(1+\beta^2)}_{2}}\int_{\Omega} \bar \xi_k |x|^{-b}\bar{D}^{1+2\beta^2}_{k}(x)dx.
\end{array}
\end{equation}

Multiplying (\ref{4.16}) by $\bar{\xi}_{k}$ and integrating over $\Omega$, we obtain that
\[\arraycolsep=1.5pt\begin{array}{lll}
&&\displaystyle \alp ^2_k\int_{\Omega} |\nabla \bar \xi_k|^2dx -\mu_{1,k}\alp ^2 _k\int_{\Omega}  |\bar\xi_k|^2dx+\alp ^2_k\int_{\Omega} V(x)|\bar\xi_k|^2 dx\\[4mm]
&=&\displaystyle \frac{a_{k}\alp^{b}_{k}\beta^{2-b}(1+2\beta^2)}{a^*}\int_{\Omega} |x|^{-b}\bar{D}^{2\beta^2}_{k}(x)|\bar\xi_k|^2dx\displaystyle\\
&&\displaystyle-\frac{2a_{k}\beta^{N+4-b}}
{\alp^{N-b}_{k}\|Q\|^{2(1+\beta^2)}_{2}}\int_{\Omega} \bar u_{2,k}\bar\xi_k dx\int_{\Omega}\bar\xi_k |x|^{-b}\bar D^{1+2\beta^2}_{k}(x)dx\\[4mm]
&\leq &\displaystyle \frac{a_{k}\alp^{b}_{k}\beta^{2-b}(1+2\beta^2)}{a^*}\int_{\Omega}|x|^{-b}\bar{D}^{2\beta^2}_{k}(x)dx
\displaystyle\\
&&\displaystyle+\frac{2a_{k}\beta^{N+4-b}}{\alp^{N-b}_{k}\|Q\|^{2(1+\beta^2)}_{2}}\int_{\Omega} \bar u_{2,k}dx\int_{\Omega}|x|^{-b}\bar{D}^{1+2\beta^2}_{k}(x)dx\\[4mm]
&\le & C\alp ^N_k\,\ \mbox{as} \ \ k\rightarrow\infty,
\end{array}\]
where we used the fact that $\bar{\xi}_{k}$ and $\tilde{u}_{i,k}(x)$ are bounded uniformly in $k$ and $\tilde{u}_{i,k}(x)$ decays exponentially as $|x|\rightarrow\infty$. Recall from (\ref{4.6}) that $\alpha_{k}^{2}\mu_{i,k}\rightarrow-\beta^2$ as $k\rightarrow\infty$, the above estimate further implies that there exists a constant $C_1>0$ such that
\begin{equation}\label{4.19}
I:=\alp ^2_k\int_{\Omega} |\nabla \bar \xi_k|^2dx+\frac{\beta ^2}{2} \int_{\Omega} |\bar\xi_k|^2dx+ \alp ^2_k\int_{\Omega} V(x)|\bar\xi_k|^2dx<C_1\alp ^N_k \quad \text{as}\ \, k\rightarrow\infty.
\end{equation}
We then conclude from \cite[Lemma 4.5]{PLC} that for any $x_0\in\Omega$, there exist constants $\delta >0$ and $C_2>0$  such that
\[
\int_{\partial B_\delta (x_0)} \Big[ \alp ^2_k |\nabla \bar \xi_k|^2+ \frac{\beta ^2}{2}  |\bar\xi_k|^2+ \alp ^2_k  V(x)|\bar\xi_k|^2\Big]dS\le C_2I\le C_1C_2\alp ^N_k\,\ \mbox{as} \,\ k\rightarrow\infty,
\]
and the claim (\ref{4.14}) is hence proved.

Under the assumptions of Theorem \ref{th4}, we complete the proof by deriving a contradiction through the following three steps.
\vskip 0.1truein
\noindent{\em  Step 1.} There exists a constant $b_0$ such that up to a subsequence if necessary, $\bar \xi_k\big(\frac{\alp _k}{\beta}x\big)\rightarrow\xi_{0}(x)$ in $C_{loc}(\mathbb{R}^{N})$ as $k\rightarrow\infty$, where
\begin{equation}\label{4.22}
\xi_{0}=b_0\big(\frac{N}{2}Q+x\cdot\nabla Q\big).
\end{equation}

To prove (\ref{4.22}), we denote
\begin{equation}\label{4.23}
\tilde{D}^{s-1}_{k}(x):=\bar{D}^{s-1}_{k}(\frac{\alpha_{k}}{\beta}x)=\int^{1}_{0}\big[t\tilde{u}_{2,k}+(1-t)\tilde{u}_{1,k}
\big]^{s-1}dt,
\end{equation}
and
\begin{equation}\label{4.20}
\xi _k(x):=\bar \xi_k\big(\frac{\alp _k}{\beta}x\big),
\end{equation}
 then $\xi_{k}$ satisfies
\begin{equation}\label{4.24}
-\Delta \xi_k(x)+C_{k}(x)\xi_k(x)=g_k(x)\quad \text{in\,\, $\Omega_{k}$},
\end{equation}
where $\lim_{k\rightarrow\infty}\Omega_{k}=\mathbb{R}^{N}$ and the coefficients
\begin{equation}\label{4.25}
C_k(x)=-\frac{a_{k}(1+2\beta^2)}{a^*|x|^b}\tilde{D}^{2\beta^2}_{k}(x)-\frac{\alp ^2_k}{\beta ^2}\mu _{1,k}+\frac{\alp ^2_k}{\beta ^2}V\big(\frac{\alp_k x}{\beta}\big),
\end{equation}
and
\begin{equation}\label{4.26}
\arraycolsep=1.5pt\begin{array}{lll}
g_k(x)&=&\displaystyle\frac{\tilde{u}_{2,k}}{\beta ^2}\frac{\alp ^2_k(\mu _{2,k}-\mu _{1,k})}{\|\tilde{u}_{2,k}-\tilde{u}_{1,k}\|_{L^\infty(\Omega_{k})}}\\[4mm]
&=&-\displaystyle \frac{2a_{k}\tilde{u}_{2,k}\beta^2}{\|Q\|_2^{2(1+\beta^2)}}\int_{\Omega_{k}} |x|^{-b}\xi_{k}(x)\tilde{D}^{1+2\beta^2}_{k}(x)dx.
\end{array}
\end{equation}

Recall that $\tilde{u}_{i,k}(x)$ is bounded uniformly in $L^{\infty}(\mathbb{R}^{N})$ ,
we then obtain that $|x|^{-b}\tilde{D}_{k}^{2\beta^{2}}(x)$ is bounded uniformly in $L_{loc}^{r}(\mathbb{R}^{N})$, where $r\in(\frac{N}{2}, \frac{N}{b})$. Since $\|\xi_{k}\|_{L^{\infty}(\Omega_{k})}\leq1$, the standard elliptic regularity (cf.\cite{GT}) then yields that $\|\xi_{k}\|_{W_{loc}^{2, r}(\Omega_{k})}\leq C$ and thus $\|\xi_{k}\|_{C_{loc}^{\alpha}(\Omega_{k})}\leq C$ for some $\alpha\in(0, 2-b)$, where the constant $C>0$ is independent of $k$. Therefore, there exist a subsequence ${a_{k}}$, still denoted by ${a_{k}}$, and a function $\xi_{0}(x)$ such that $\xi_{k}(x)\rightarrow\xi_{0}(x)$ in $C_{loc}(\mathbb{R}^{N})$ as $k\rightarrow\infty$. Applying Lemma \ref{lemma3.2} and (\ref{4.6}), direct calculations yield from (\ref{4.25}) and (\ref{4.26}) that
\[
C_k(x)\to 1-(1+2\beta^2)|x|^{-b}Q^{2\beta^2}(x)\ \ \text{in\ \ $C_{loc}(\R^N)$}\,\ \mbox{as} \,\ k\rightarrow\infty,
\]
and
\[
g_k(x)\to -\frac{2Q(x)\beta^2}{\|Q\|^{2}_{2}}\int_{\R^{N}}|x|^{-b}\xi _0Q^{1+2\beta^2}(x)dx\ \ \text{in\ \ $C_{loc}(\R^N)$}\,\ \mbox{as} \,\ k\rightarrow\infty.
\]
Hence, we derive from (\ref{4.24}) that $\xi_{0}$ solves
\begin{equation}\label{4.27}
\mathcal{L}\xi_0=-\Delta \xi_0+\xi_0-(1+2\beta^2)|x|^{-b}Q^{2\beta^2}\xi_0
=-\frac{2Q(x)\beta^2}{\|Q\|^{2}_{2}}\int_{\R^{N}}|x|^{-b}\xi _0Q^{1+2\beta^2}\ \ \mbox{in} \ \ \R^N.
\end{equation}
On the other hand, recall from \cite[Lemma C.1]{LZ} that $\mathcal{L}\big(\frac{N}{2}Q+x\cdot\nabla Q\big)=-2Q$, then we conclude from the non-degeneracy of $\mathcal{L}$ in (\ref{4.1}) that (\ref{4.22}) holds for some constant $b_0\in\mathbb{R}$.
\vskip 0.1truein
\noindent{\em  Step 2.} The constant $b_0=0$ in (\ref{4.22}).
\vskip 0.1truein
Similar to the argument of (\ref{3.19}), we can deduce that $\|\tilde{u}_{i,k}(x)\|_{W_{loc}^{2,r}(\Omega_{k})}\leq C$, where $r\in(\frac{N}{2},\frac{N}{b})$ and $C>0$ is independent of $k$.
Since $0<b<min\{2,\frac{N}{2}\}$, we further obtain that $\|\tilde{u}_{i,k}(x)\|_{H_{loc}^{2}(\Omega_{k})}\leq C$.
Applying the Cauchy-Schwarz inequality and the Young's inequality, it follows that
\begin{equation*}
\arraycolsep=1.5pt
\begin{array}{lll}
&&\displaystyle\alp _k^2\int_{B_{\delta}(0)} (x\cdot\nabla\bar{u}_{i,k}(x))\Delta\bar{u}_{i,k}(x)dx\\
&=&\displaystyle\alp _k^2\int_{B_{\delta}(0)} (x\cdot\nabla\tilde{u}_{i,k}(\frac{\beta}{\alp _k}x))\Delta\tilde{u}_{i,k}(\frac{\beta}{\alp _k}x)dx\\
&=&\displaystyle\frac{\alp _k^N}{\beta^{N-2}}\int_{B_{\frac{\beta\delta}{\alp _k}}(0)} (x\cdot\nabla\tilde{u}_{i,k}(x))\Delta\tilde{u}_{i,k}(x)dx\\
&\leq&\displaystyle\frac{\delta\alp _k^{N-1}}{\beta^{N-3}}\int_{B_{\frac{\beta\delta}{\alp _k}}(0)} |\nabla\tilde{u}_{i,k}(x)||\Delta\tilde{u}_{i,k}(x)|dx\\
&\leq&\displaystyle\frac{\delta\alp _k^{N-1}}{2\beta^{N-3}}\Big(\int_{B_{\frac{\beta\delta}{\alp _k}}(0)}|\nabla\tilde{u}_{i,k}(x)|^{2}dx+\int_{B_{\frac{\beta\delta}{\alp _k}}(0)}|\Delta\tilde{u}_{i,k}(x)|^{2}dx\Big)\\
&\leq&\displaystyle\frac{\delta\alp _k^{N-1}}{2\beta^{N-3}}\|\tilde{u}_{i,k}(x)\|_{H_{B_{\frac{\beta\delta}{\alp _k}}(0)}^{2}}^{2}\leq C,
\end{array}
\end{equation*}
where $\delta>0$ is sufficiently small such that $B_\delta(0)\subset\Omega$. Therefore, we can use the integration by parts to get that
\begin{equation}\label{4.28}
\arraycolsep=1.5pt\begin{array}{lll}
&&-\displaystyle\alp _k^2 \int_{B_{\delta}(0)}(x\cdot \nabla \bar u_{i,k}) \Delta \bar u_{i,k}dx \\[4mm]
&=&-\alp _k^2\displaystyle  \int_{\partial B_{\delta}(0)} \frac{\partial\bar u_{i,k} }{\partial \nu }(x\cdot \nabla \bar u_{i,k})dS
+\displaystyle\alp _k^2 \int_{B_{\delta}(0)} \nabla \bar u_{i,k}\cdot\nabla (x\cdot \nabla \bar u_{i,k})dx \\[4mm]
&=&\alp _k^2\displaystyle  \int_{\partial B_{\delta}(0)} \Big[-\frac{\partial\bar u_{i,k} }{\partial \nu }(x\cdot \nabla \bar u_{i,k})+\frac 12(x\cdot \nu) |\nabla \bar u_{i,k}|^2\Big]dS
+\frac{(2-N)\alp_k^2}{2}\int_{B_{\delta}(0)}|\nabla \bar u_{i,k}|^2dx\\[4mm]
&=&\alp _k^2\displaystyle  \int_{\partial B_{\delta}(0)} \Big[-\frac{\partial\bar u_{i,k} }{\partial \nu }(x\cdot \nabla \bar u_{i,k})+\frac 12(x\cdot \nu) |\nabla \bar u_{i,k}|^2+\frac{(2-N)}{4}(\nabla \bar u_{i,k}^2\cdot\nu) \Big]dS\\[4mm]
&&-\displaystyle\frac{(2-N)}{2} \int_{B_{\delta}(0)}\Big[\alp _k^2V(x)\bar u^{2}_{i,k}-\alp _k^2\mu_{i,k}\bar u_{i,k}^{2}- \frac{a_{k}\alp _k^b\beta^{2-b}\bar u_{i,k}^{2(1+\beta^2)}}{a^*|x|^b}\Big]dx,
\end{array}
\end{equation}
where the last equality has used the following fact
\begin{equation}\label{4.29}
\arraycolsep=1.5pt
\begin{array}{lll}
\displaystyle\alp _k^2\int_{B_{\delta}(0)} |\nabla \bar u_{i,k}|^2
&=&\displaystyle\frac{\alp _k^2}{2}\int_{\partial B_{\delta}(0)}(\nabla \bar u_{i,k}^2\cdot\nu) dS\\
&&\displaystyle-\int_{B_{\delta}(0)}\Big[\alp _k^2V(x)\bar u^{2}_{i,k}-\alp _k^2\mu_{i,k}\bar u_{i,k}^{2}- \frac{a_{k}\alp _k^b\beta^{2-b}\bar u_{i,k}^{2(1+\beta^2)}}{a^*|x|^b}\Big]dx.
\end{array}
\end{equation}

On the other hand, multiplying (\ref{4.12}) by  $ (x\cdot \nabla  \bar{u}_{i,k})$ and integrating over $B_\delta (0)$, where $\delta >0$ is small as before, we deduce that
\begin{equation}\label{4.30}
\arraycolsep=1.5pt\begin{array}{lll}
&&-\displaystyle\alp _k^2 \int_{B_{\delta}(0)} \big(x\cdot \nabla \bar u_{i,k}\big) \Delta \bar u_{i,k}dx \\[4mm]
&=& \displaystyle\alp _k^2 \int_{B_{\delta}(0)} \big[\mu_{i,k}-V(x)\big] \bar u_{i,k}\big(x\cdot \nabla \bar u_{i,k}\big)dx\\
&&+\displaystyle \frac{a_{k}\alp^{b}_{k}\beta^{2-b} }{a^*} \int_{B_{\delta}(0)} \bar u_{i,k}^{1+2\beta^2}|x|^{-b}\big (x\cdot \nabla \bar u_{i,k}\big)dx\\[4mm]
&=& -\displaystyle\frac{\alp _k^2}{2} \int_{B_{\delta}(0)}\bar u_{i,k}^2\Big\{N\big[\mu_{i,k}-V(x)\big]-[x\cdot \nabla V(x)]\Big\}dx\\[4mm]
&&+\displaystyle\frac{\alp _k^2}{2}  \int_{\partial B_{\delta}(0)}\bar u_{i,k}^2\big[\mu_{i,k}-V(x)\big](x\cdot\nu) dS\\[4mm]
&&+\displaystyle \frac{a_{k}\alp^{b}_{k}\beta^{2-b} }{2a^*(1+\beta^2)} \Big[\displaystyle \int_{\partial B_{\delta}(0)} \bar{u}_{i,k}^{2(1+\beta^2)}\frac{(x\cdot\nu)}{|x|^{b}} dS-\int_{B_{\delta}(0)} \bar u_{i,k}^{2(1+\beta^2)}\frac{N-b}{|x|^{b}}dx\Big].\\[4mm]
\end{array}
\end{equation}

Substituting (\ref{4.28}) into (\ref{4.30}), one can derive that
\begin{equation}\label{4.31}
\arraycolsep=1.5pt\begin{array}{lll}
&&\alp _k^2\displaystyle  \int_{\partial B_{\delta}(0)} \Big[-\frac{\partial\bar u_{i,k} }{\partial \nu }(x\cdot \nabla \bar u_{i,k})+\frac 12(x\cdot \nu) |\nabla \bar u_{i,k}|^2+\frac{(2-N)}{4}(\nabla \bar u_{i,k}^2\cdot\nu) \Big]dS\\[4mm]
&=&\displaystyle \int_{B_{\delta}(0)}\Big\{\alp _k^2\big[-\mu_{i,k}+V(x)+\frac{1}{2}(x\cdot\nabla V)\big] \bar u_{i,k}^2-\frac{a_{k}\alp^{b}_{k}\beta^{4-b} }{a^*(1+\beta^2)}\bar u_{i,k}^{2(1+\beta^2)}|x|^{-b}\Big\}dx+I_i,
\end{array}
\end{equation}
where $I_i$ satisfies
\begin{equation}\label{4.32}
\arraycolsep=1.5pt\begin{array}{lll}
I_i&=&\displaystyle\frac{\alp _k^2}{2} \int_{\partial B_{\delta}(0)}\bar u_{i,k}^2\big[\mu_{i,k}-V(x)\big](x\cdot\nu) dS\\
&&+\displaystyle \frac{a_{k}\alp^{b}_{k}\beta^{2-b} }{2a^*(1+\beta^2)} \int_{\partial B_{\delta}(0)} \bar u_{i,k}^{2(1+\beta^2)}\frac{(x\cdot\nu)}{|x|^{b}}dS.
\end{array}
\end{equation}
Note from (\ref{4.2}) that
\[\arraycolsep=1.5pt\begin{array}{lll}
&&\quad\displaystyle \mu_{i,k}\displaystyle\alp _k^2 \int_{\Omega} \bar u_{i,k}^2dx+\displaystyle \frac{a_{k}\alp^{b}_{k}\beta^{4-b} }{a^*(1+\beta^2)} \int_{\Omega}  \bar u_{i,k}^{2(1+\beta^2)}|x|^{-b}dx=\displaystyle\frac{\|Q\|^{2}_{2}\alp^{N+2}_{k}}{\beta^{N}}e(a_{k}).
\end{array}\]
We then deduce from above that
\[\arraycolsep=1.5pt\begin{array}{lll}
&&\quad\displaystyle-\alp^{2}_{k}\int_{B_{\delta}(0)} \Big[V(x)+\frac{1}{2}[x\cdot \nabla V(x)]\Big]\bar u_{i,k}^{2}dx+\frac{\|Q\|^{2}_{2}\alp^{N+2}_{k}}{\beta^{N}}e(a_{k})\\[4mm]
&&=I_i+\alp _k^2\displaystyle \int_{\partial B_{\delta}(0)} \frac{\partial\bar u_{i,k} }{\partial \nu }(x\cdot \nabla \bar u_{i,k})dS
-\displaystyle \frac{\alp _k^2}{2}\int_{\partial B_{\delta}(0)} (x\cdot \nu )|\nabla \bar u_{i,k}|^2dS\\[4mm]
&&\quad-\displaystyle\frac{2-N}{4}\alp _k^2\int_{\partial B_{\delta}(0)} \nabla {\bar u_{i,k}}^{2}\cdot\nu dS+
\mu_{i,k}\alp _k^2 \int _{\Omega\backslash B_\delta (0)} \bar u_{i,k}^2dx\\[4mm]
&&\quad+\displaystyle \frac{a_{k}\alp^{b}_{k}\beta^{4-b} }{a^*(1+\beta^2)} \int _{\Omega\backslash B_\delta (0)} \bar u_{i,k}^{2(1+\beta^2)}|x|^{-b}dx,
\end{array}\]
which further implies that
\begin{equation}\label{4.33}
\displaystyle-\alp^{2}_{k}\int_{B_{\delta}(0)} \Big[V(x)+\frac{1}{2}[x\cdot \nabla V(x)]\Big](\bar u_{1,k}+\bar u_{2,k})\bar{\xi}_k=T_k,
\end{equation}
here $T_k$ satisfies
\begin{equation}\label{4.34}
\arraycolsep=1.5pt\begin{array}{lll}
T_k&=&\displaystyle\frac{I_2-I_1}{\|\bar u_{2,k}-\bar u_{1,k}\|_{L^\infty(\Omega)}}-\displaystyle\frac{\alp _k^2}{2}\int_{\partial B_{\delta}(0)} (x\cdot \nu )\big[\big(\nabla \bar u_{2,k}+\nabla \bar u_{1,k}\big)\cdot\nabla \bar \xi_k\big]dS\\[4mm]
&&-\displaystyle \frac{2-N}{4}\alp _k^2\Big[\int_{\partial B_{\delta}(0)}\big[(\nabla \bar u_{2,k}
+\nabla \bar u_{1,k}\big)\cdot\nu] \bar \xi_k dS\\
&&+\displaystyle\int_{\partial B_{\delta}(0)}\big(\bar u_{2,k}+ \bar u_{1,k}\big) (\nabla\bar \xi_k\cdot\nu)dS\Big]\\[4mm]

&&+\displaystyle  \alp _k^2 \int_{\partial B_{\delta}(0)} \Big[(x\cdot \nabla \bar u_{2,k})\big(\nu \cdot \nabla\bar \xi_k \big)+\big(\nu \cdot \nabla \bar u_{1,k}\big)(x \cdot \nabla \bar \xi_k)\Big]dS
\\[4mm]
&&+\mu_{2,k}\displaystyle\alp _k^2 \int _{\Omega\backslash B_\delta (0)} \big(\bar u_{2,k}+\bar u_{1,k}\big)\bar \xi_k dx+\frac{(\mu_{2,k}-\mu_{1,k})\alp_k^2}{\|\bar u_{2,k}-\bar u_{1,k}\|_{L^\infty(\Omega)}}\int _{\Omega\backslash B_\delta (0)} \bar u_{1,k}^2dx\\[4mm]
&&+\displaystyle\frac{2a_{k}\alp^{b}_{k}\beta^{4-b}}{a^*} \int _{\Omega\backslash B_\delta (0)} \bar{D}^{1+2\beta^2}_{k}|x|^{-b}\bar \xi_k dx,\\[4mm]
\end{array}
\end{equation}
where
\begin{equation}\label{4.35}
\arraycolsep=1.5pt\begin{array}{lll}
\displaystyle\frac{I_2-I_1}{\|\bar u_{2,k}-\bar u_{1,k}\|_{L^\infty(\Omega)}}
&=&\displaystyle \frac{a_{k}\alp^{b}_{k}\beta^{2-b}}{a^*} \int_{\partial B_{\delta}(0)} (x\cdot\nu)|x|^{-b} \bar{D}^{1+2\beta^2}_{k} \bar \xi_k  dS\\[4mm]
&&-\displaystyle\frac{\alp _k^2}{2} \int_{\partial B_{\delta}(0)} \big(\bar u_{2,k}+\bar u_{1,k}\big)\bar \xi_k V(x)(x\cdot\nu) dS\\[4mm]
&&+\displaystyle\frac{\mu_{2,k}\alp _k^2}{2} \int_{\partial B_{\delta}(0)}   \big(\bar u_{2,k}+\bar u_{1,k}\big)\bar \xi_k (x\cdot\nu) dS\\[4mm]
&&+\displaystyle\frac{\big(\mu_{2,k}-\mu_{1,k}\big)\alp _k^2}{2\|\bar u_{2,k}-\bar u_{1,k}\|_{L^\infty(\Omega)}} \int_{\partial B_{\delta}(0)}   \bar u_{1,k}^2 (x\cdot\nu) dS.
\end{array}
\end{equation}

We next estimate the right hand side of (\ref{4.34}). We first consider
$\displaystyle\frac{I_2-I_1}{\|\bar u_{2,k}-\bar u_{1,k}\|_{L^\infty(\Omega)}}$. Applying the H\"{o}lder inequality, it follows from the exponential decay of $\tilde u_{i,k}(x)$ and (\ref{4.14}) that
\begin{equation}\label{4.36}
\arraycolsep=1.5pt\begin{array}{lll}
&&\displaystyle \frac{a_{k}\alp^{b}_{k}\beta^{2-b}}{a^*} \int_{\partial B_{\delta}(0)} (x\cdot\nu)|x|^{-b} \bar{D}^{1+2\beta^2}_{k} \bar \xi_k  dS\\
&&-\displaystyle\frac{\alp _k^2}{2} \int_{\partial B_{\delta}(0)} \big(\bar u_{2,k}+\bar u_{1,k}\big)\bar \xi_k V(x)(x\cdot\nu) dS\\
&&+\displaystyle\frac{\mu_{2,k}\alp _k^2}{2} \int_{\partial B_{\delta}(0)}   \big(\bar u_{2,k}+\bar u_{1,k}\big)\bar \xi_k (x\cdot\nu) dS=o(e^{-\frac{C\delta}{\alp_k}})\ \ \text{as} \ \ k\rightarrow\infty.
\end{array}
\end{equation}
Moreover, we deduce from (\ref{4.18}) that
\begin{equation}\label{4.37}
\displaystyle\frac{\big|\mu_{2,k}-\mu_{1,k}\big|\alp _k^2}{\|\bar u_{2,k}-\bar u_{1,k}\|_{L^\infty(\Omega)}}
\le  \displaystyle\frac{2a_{k}\beta^{N+4-b}}{\alp^{N-b}_{k}\|Q\|^{2(1+\beta^2)}_{2}}\int_{\Omega} |x|^{-b}\bar{D}^{1+2\beta^2}_{k}|\bar \xi_k |
\le  C\ \ \text{as}\ \ k\rightarrow\infty,
\end{equation}
which combined with (\ref{4.36}) yield that
\begin{equation}\label{4.38}
\displaystyle\frac{I_2-I_1}{\|\bar u_{2,k}-\bar u_{1,k}\|_{L^\infty(\Omega)}}=o(e^{-\frac{C\delta}{\alp_k}})\ \ \text{as}\ \ k\rightarrow\infty.
\end{equation}

In addition, using the H\"{o}lder inequality, we derive from (\ref{4.9}) and (\ref{4.14}) that
\begin{equation*}
\arraycolsep=1.5pt\begin{split}
&\alp_{k}^2\int_{\partial B_{ \delta}(0)}\big|(x\cdot\nabla \bar{u}_{2,k})(\nu\cdot\nabla\bar{\xi}_{k})\big|dS\\
\leq& \alp_{k}\big(\int_{\partial B_{ \delta}(0)}\big|\nabla \bar{u}_{2,k}\big|^2dS\big)^{\frac{1}{2}}
\big(\alp_{k}^2\int_{\partial B_{ \delta}(0)}\big|\nabla \bar{\xi}_{k}\big|^2dS\big)^{\frac{1}{2}}\\
\leq& C\alp_{k}^{\frac{N}{2}+1}e^{-\frac{C\delta}{\alp_k}}\ \ \text{as} \ \ k\rightarrow\infty,
\end{split}
\end{equation*}
where $C>0$ is independent of $k$. Similarly, one can get that
\begin{equation*}
\arraycolsep=1.5pt\begin{split}
&\alp_{k}^2\int_{\partial B_{ \delta}(0)}\big|(\nu\cdot\nabla \bar{u}_{1,k})(x\cdot\nabla\bar{\xi}_{k})\big|dS\\
+&\alp_{k}^2\int_{\partial B_{\delta}(0)} \big|(x\cdot \nu )\big[\big(\nabla \bar u_{2,k}+\nabla \bar u_{1,k}\big)\cdot\nabla \bar \xi_k\big]\big|dS\\
+&\alp_{k}^2\int_{\partial B_{ \delta}(0)}\big|( \bar{u}_{1,k}+ \bar{u}_{2,k})(\nu\cdot\nabla\bar{\xi}_{k})\big|dS\\
+&\alp_{k}^2\int_{\partial B_{\delta}(0)}\big|\nu\cdot( \nabla\bar{u}_{2,k}+ \nabla\bar{u}_{1,k})\bar{\xi}_{k}\big|dS \leq C\alp_{k}^{\frac{N}{2}+1}e^{-\frac{C\delta}{\alp_k}}
\ \ \text{as} \ \ k\rightarrow\infty.
\end{split}
\end{equation*}
On the other hand, since $|\bar{\xi}_{k}|$ and $|\mu_{i,k}\alpha_{k}^{2}|$ are bounded uniformly in $k$, we deduce that
\begin{equation*}
\arraycolsep=1.5pt\begin{split}
&\Big|\mu_{2,k}\displaystyle\alp _k^2 \int _{\Omega\backslash B_\delta (0)} \big(\bar u_{2,k}+\bar u_{1,k}\big)\bar \xi_kdx+\frac{(\mu_{2,k}-\mu_{1,k})\alp_k^2}{\|\bar u_{2,k}-\bar u_{1,k}\|_{L^\infty(\Omega)}}\int _{\Omega\backslash B_\delta (0)} \bar u_{1,k}^2dx\\
&+\displaystyle \frac{2a_{k}\alp^{b}_{k}\beta^{4-b} }{a^*} \int _{\Omega\backslash B_\delta (0)} \bar{D}^{1+2\beta^2}_{k}|x|^{-b}\bar \xi_kdx\Big|=o(e^{-\frac{C\delta}{\alp_k}})\ \ \text{as} \ \ k\rightarrow\infty.
\end{split}
\end{equation*}
We finally conclude from above  that
\[
T_k=o(e^{-\frac{C\delta}{\alp_k}}) \,\ \mbox{as} \,\ k\rightarrow\infty,
\]
and we obtain from \eqref{4.33} that
\begin{equation}\arraycolsep=1.5pt\begin{array}{lll}
\displaystyle-\alp^{2}_{k}\int_{B_{\delta}(0)} \Big[V(x)+\frac{1}{2}[x\cdot \nabla V(x)]\Big](\bar u_{1,k}+\bar u_{2,k})\bar{\xi}_k=o(e^{-\frac{C\delta}{\alp_k}})\ \ \text{as}\ \ k\rightarrow\infty.
\end{array} \label{4.39}
\end{equation}

Furthermore, under the assumptions of Theorem \ref{th4}, one can check that
\begin{equation*}
V(x)=[C_{0}+o(1)]|x|^{p},\ \ \frac{\partial V(x)}{\partial x_{j}}=[C_{0}+o(1)]\frac{\partial|x|^{p}}{\partial x_{j}}\ \ \text{as}\ \ |x|\rightarrow 0,
\end{equation*}
where $j=1,\cdots,N,$ and $C_{0}$ is a positive constant satisfying $\lim_{x\rightarrow 0}\frac{V(x)}{|x|^{p}}=C_{0}$.
Therefore, for $x\in B_{\delta}(0)$, where $\delta>0$ is small enough, we have
\[
V(x)+\frac{1}{2}[x\cdot \nabla V(x)]=\big[C_{0}+o(1)\big]\big[|x|^{p}+\frac{1}{2}(x\cdot \nabla |x|^{p})\big]=[C_{0}+o(1)]\frac{2+p}{2}|x|^{p},
\]
where the last equality follows from the fact that $x\cdot \nabla |x|^{p}=p|x|^{p}$. We then derive that
\begin{equation}\label{4.40}
\arraycolsep=1.5pt\begin{split}
o(e^{-\frac{C\delta}{\alp_k}})=&\displaystyle-\alp^{2}_{k}\int_{B_{\delta}(0)} \Big[V(x)+\frac{1}{2}\big(x\cdot \nabla V(x)\big)\Big](\bar u_{1,k}+\bar u_{2,k})\bar{\xi}_kdx\\
=&-\alp^{2}_{k}\int_{B_{\delta}(0)}[C_{0}+o(1)]\frac{2+p}{2}|x|^{p} \big(\bar u_{1,k}(x)+\bar u_{2,k}(x)\big)\bar{\xi}_kdx\\
=&\displaystyle-\frac{\alp^{2+N+p}_{k}}{\beta^{N+p}}\int_{B_{\frac{\beta\delta}{\alp_{k}}}(0)} [C_{0}+o(1)]\frac{2+p}{2} |y|^{p}\big(\tilde u_{1,k}(y)+\tilde u_{2,k}(y)\big)\xi_{k}(y)dy\\
&\ \ \text{as}\ \ k\rightarrow\infty.
\end{split}
\end{equation}
 Note that $\tilde u_{i,k}(x)=\bar u_{i,k}\big(\frac{\alp_k}{\beta}x\big)\to Q(x)$ uniformly in $\mathbb{R}^{N}$ as $k\rightarrow\infty$, we deduce from (\ref{4.22}) and (\ref{4.40}) that
\begin{equation*}
\arraycolsep=1.5pt\begin{split}
0&=\int_{\mathbb{R}^{N}}|x|^{p}Q(x)\xi_0(x)dx\\
&=b_0\int_{\mathbb{R}^{N}}|x|^{p}Q(x)\big(\frac{N}{2}Q+x\cdot\nabla Q\big)dx\\
&=\frac{Nb_0}{2}\int_{\mathbb{R}^{N}}|x|^{p}Q^2(x)dx
+\frac{b_0}{2}\int_{\mathbb{R}^{N}}|x|^{p}\big(x\cdot \nabla Q^2(x)\big)dx\\
&=-\frac{pb_0}{2}\int_{\mathbb{R}^{N}}|x|^{p}Q^2(x)dx,
\end{split}
\end{equation*}
which implies that $b_0=0$, i.e., $\xi_0\equiv0$.
\vskip 0.1truein
\noindent{\em  Step 3.}
$\xi_0\equiv 0$ cannot occur.
\vskip 0.1truein
Let $y_k$ be a point satisfying $|\xi_k(y_k)|=1$.
 Since $\tilde u_{i,k}(x)$ and $Q(x)$ decay exponentially as $k\rightarrow\infty$, applying the maximum principle to (\ref{4.24}) yields that $|y_k|\le C$ uniformly in $k$.  Therefore, we have $\xi_{k}\rightarrow\xi_{0}\not\equiv 0$ in $C_{loc}(\mathbb{R}^{N})$,
 which however contradicts the fact that $\xi_0 \equiv 0$. This completes the proof of Theorem \ref{th4}.
\qed
\vspace {.1cm}
\begin{appendix}
\section{Appendix. }
In this section, we shall prove Proposition \ref{2.1} on the Gagliardo-Nirenberg type inequality defined in any open bounded domain.

\vskip 0.1truein
\noindent\textbf{Proof of Proposition \ref{2.1}.} For any given open bounded domain $\Omega\subset\mathbb{R}^{N}$, stimulated by \cite{NTV}, we define
\begin{equation}\label{A1}
\frac{a^{*}(\Omega)}{1+\beta^{2}}:=\inf_{\{u\neq0, u\in H_{0}^{1}(\Omega)\}}\Upsilon(u),
\end{equation}
where
\begin{equation}\label{A2}
\Upsilon(u):=\frac{\displaystyle\int_{\Omega}|\nabla u(x)|^{2}dx\Big(\int_{\Omega}|u(x)|^{2}dx\Big)^{\beta^{2}}}{\displaystyle\int_{\Omega}|x|^{-b}|u(x)|^{2+2\beta^{2}}dx}.
\end{equation}
To complete the proof, we only need to prove that $a^{*}(\Omega)=a^{*}:=\|Q\|_{L^{2}(\mathbb{R}^{N})}^{2\beta^{2}}$ and the infimum (\ref{A1}) is not attained.

We first recall from \cite{AI} the following Gagliardo-Nirenberg inequality:
\begin{equation}\label{A3}
\frac{a^{*}}{1+\beta^{2}}\int_{\mathbb{R}^{N}}|x|^{-b}|u(x)|^{2+2\beta^{2}}dx\leq\int_{\mathbb{R}^{N}}|\nabla u(x)|^{2}dx\Big(\int_{\mathbb{R}^{N}}|u(x)|^{2}dx\Big)^{\beta^{2}}\ \text{ in}\,\ \ H^{1}(\mathbb{R}^N).
\end{equation}
It is easy to infer that $a^{*}\leq a^{*}(\Omega)$, since $H_{0}^{1}(\Omega)\subset H^{1}(\mathbb{R}^{N})$. We next prove that $a^{*}\geq a^{*}(\Omega)$. Since $0\in\Omega$ is an open bounded domain in $\mathbb{R}^{N}$, there exists an open ball $B_{2R}(0)\subset\Omega$ centered at an inner point $0$, where $R>0$ is small enough.
Choose a nonnegative cutoff function $\varphi\in C_{0}^{\infty}(\mathbb{R}^{N})$ such that $\varphi(x)=1$ for $|x|\leq R$ and $\varphi(x)=0$ for $|x|\geq2R$. For $\tau>0$ and $A_{\tau}>0$, set
\begin{equation}\label{A4}
\Phi_{\tau}(x)=A_{\tau}\frac{\tau^{\frac{N}{2}}}{\|Q\|_{2}}\varphi(x)Q(\tau x),\ \ x\in\Omega,
\end{equation}
where $Q=Q(|x|)$ is the unique positive radial solution of (\ref{1.5}) and $A_{\tau}>0$ is chosen such that $\int_{\Omega}|\Phi_{\tau}(x)|^{2}dx=1$. One can check that
\begin{equation}\label{A5}
\frac{a^{*}(\Omega)}{1+\beta^{2}}\leq \Upsilon(\Phi_{\tau})\rightarrow \frac{a^{*}}{1+\beta^{2}}\ \text{ as}\,\ \ \tau\rightarrow\infty,
\end{equation}
which implies that $a^{*}\geq a^{*}(\Omega)$. Therefore, we conclude that $a^{*}=a^{*}(\Omega)$.

Finally, we prove that the infimum (\ref{A1}) is not attained. On the contrary, suppose $\Upsilon(u)$ attains its infimum at some nonzero $u_{0}\in H_{0}^{1}(\Omega)$, it then implies that the equality (\ref{A3}) is attained at $\bar{u}(x)\in H^{1}(\mathbb{R}^{N})$, where $\bar{u}(x)=u_{0}(x)$ for $x\in\Omega$ and $\bar{u}(x)\equiv0$ for $x\in\mathbb{R}^{N}\backslash \Omega$. However, this contradicts the fact that the equality of (\ref{A3}) holds if and only if $u(x)=mn^{\frac{N}{2}}Q(nx)$ ($m, n\neq 0$ are arbitrary). Therefore, the proof of Proposition \ref{2.1} is complete.\qed

\end{appendix}

\vspace {.5cm}

\noindent {\bf Acknowledgements:}  The authors thank Professor Yujin Guo very much for his fruitful discussions on the present paper.

\vspace {.5cm}
\noindent {\bf Data Availability}

Data sharing is not applicable to this article as no new data were created or analyzed in this study.


\begin{thebibliography}{45}


\bibitem{GPA}  G. Agrawal, {\em Nonlinear fiber optics}, Lecture Notes in Physics, \textbf{18} (2005).

\bibitem{AA}  A. Ardila, V. Dinh, {\em Some qualitative studies of the focusing inhomogeneous Gross-Pitaevskii equation}, Z. Angew. Math. Phys. \textbf{71} (2020), 79.



\bibitem{BP}  G. Baym, C. Pethick, {\em  Ground state properties of magnetically trapped Bose Einstein
    condensate rubidium gas}, Phys. Rev. Lett. \textbf{76} (1996), 6-9.

\bibitem{PLC}  D. Cao, S. Li and P. Luo, {\em Uniqueness of positive bound states with multi-bump for nonlinear Schr\"{o}dinger equations}, Calc. Var. Paratial Differ. Equ. \textbf{54} (2015), 4037-4063.

\bibitem{YPC} D. Cao, S. Peng and S. Yan, Singularly Perturbed Methods for Nonlinear Elliptic Problems, Cambridge University Press,  Cambridge, 2021.

\bibitem{DGL1}  Y. Deng,  Y. Guo and L. Lu, {\em On the collapse and concentration of Bose-Einstein condensates with inhomogeneous attractive interactions}, Calc. Var. Partial Differ. Equ. {\bf 54} (2015), 99-118.

\bibitem{DGL2}  Y. Deng,  Y. Guo and L. Lu, {\em Threshold behavior and uniqueness of ground states for mass critical inhomogeneous Schrdinger equations}, J. Math. Phys. {\bf 59} (2018), 011503.





\bibitem{F}  F. Genoud, Th\'{e}orie de bifurcation et de stabilit\'{e} pour une \'{e}quation de Schr\"{o}dinger avec une non-lin\'{e}arit\'{e} compacte. Ph. D thesis. EPFL, 2008.

\bibitem{G}  F. Genoud, {\em A uniqueness result for $\Delta u-\lambda u+V(x)u^{p}=0$ on $\mathbb{R}^{2}$}, Adv. Nonlinear Stud. \textbf{11} (2011), 483-491.

\bibitem{AI}  F. Genoud, {\em An inhomogeneous, $L^{2}$-critical, nonlinear schr\"{o}dinger equation}, Z. Anal. Anwend. \textbf{31} (2012), 283-290.

\bibitem{WLG}  Y. Guo, C. Lin and J. Wei, {\em Local uniqueness and refined spike profiles of ground states for two-dimensional attractive Bose-Einstein condensates}, SIAM J. Math. Anal. \textbf{49} (2017), 3671-3715.


\bibitem{GL}  Y. Guo, Y. Luo and Q. Zhang, {\em Minimizers of mass critical Hartree energy functionals in bounded domains}, J. Differential Equations \textbf{265} (2018), 5177-5211.

\bibitem{SC}  F. Genoud, C. Stuart, {\em Schr\"{o}dinger equations with a spatially decaying nonlinearity: existence and stability of standing waves}, Discrete Contin. Dyn. Syst. \textbf{21} (2008), 137-186.

\bibitem{GS}  Y. Guo and R. Seiringer, {\em On the mass concentration for Bose-Einstein condensates with attractive interactions}, Lett. Math. Phys. \textbf{104} (2014), 141-156.

\bibitem{GT} D. Gilbarg, N. S. Trudinger, Elliptic Partial Differential Equations of Second Order, Springer, 1997.

\bibitem{GWZ}  Y. Guo, Z. Wang, X. Zeng and H. Zhou, {\em Properties of ground states of
attractive Gross-Pitaevskii equations with multi-well potentials}, Nonlinearity {\bf 31}
(2018), 957-979.

\bibitem{GZZ}   Y. Guo, X. Zeng and H. Zhou, {\em Energy estimates and symmetry breaking
in attactive Bose-Einstein condensates with ring-shaped potentials}, Ann. Inst. H.
Poincar\'{e} Anal. Non Lin\'{e}aire \textbf{33} (2016), 809-828.

\bibitem{HL} Q. Han, F. Lin, Elliptic Partial Differential Equations, Courant Lecture Note
in Math. 1, Courant Institute of Mathematical Science/AMS, New York, 2011.


\bibitem{LE}  E. H. Lieb and M. Loss, Analysis, $2^{nd}$ edition, Gradudate Studies in Mathematics, Amer. Math. Soc., Providence, 2001.

\bibitem{LL}  Y. Li, Y. Luo, {\em Existence and uniqueness of ground states for attractive  Bose-Einstein condensates in box-shaped traps}, J. Math. Phys. \textbf{62} (2021), 031513.

\bibitem{LT}  C. Liu, V. K. Tripathi, {\em Laser guiding in an axially nonuniform plasma channel}, Phys. Plasmas \textbf{1} (1994), 3100-3103.





\bibitem{LZ}  Y. Luo, S. Zhang, {\em Concentration behavior of ground states for $L^{2}$-critical schr\"{o}dinger equation with a spatially decaying nonlinearity}, Comm. Pure Appl. Anal. \textbf{21} (2022), 1481-1504.


\bibitem{ZL}  Y. Luo, X. Zhu, {\em Mass concentration behavior of Bose-Einstein condensates with attractive interactions in bounded domains}, Appl. Anal. \textbf{99} (2020), 2414-2427.

\bibitem{NTV}  B. Noris, H. Tavares and G. Verzini, {\em Existence and orbital stability of the ground states with prescribed mass for the $L^{2}$-critical and supercritical NLS on bounded domains}, Anal PDE. \textbf{7} (2015), 1807-1838.




\bibitem{TJ}  J. Toland, {\em Uniqueness of positive solutions of some semilinear Sturm-Liouville problems on the half line}, Proc. Roy. Soc. Edinburgh Sect. A \textbf{97}  (1984), 259-263.


\end{thebibliography}
\end{document}